\documentclass[english]{elsarticle}
\usepackage[T1]{fontenc}
\usepackage[latin9]{inputenc}
\usepackage{geometry}
\usepackage{verbatim}
\usepackage{amsmath}
\usepackage{float}
\usepackage{amsbsy}
\usepackage{amssymb}
\usepackage{stackrel}
\usepackage{graphicx}
\usepackage{esint}
\usepackage{bm}
\usepackage{algpseudocode}
\usepackage{algorithm}
\usepackage{courier}

\usepackage{amsfonts} 

\usepackage{textcomp}
\usepackage{tabularx}
\usepackage{lineno}

\usepackage{todonotes}

\usepackage{enumitem}

\usepackage[table,xcdraw]{xcolor}
\usepackage{adjustbox,lipsum}
\usepackage{subcaption}

\usepackage{xcolor}

\newcommand{\zero}{\bm{0}}

\newcommand{\R}{\mathbb{R}}
\newcommand{\N}{\mathbb{N}}

\newcommand{\q}{\bm{q}}
\newcommand{\qd}{\dot{\bm{q}}}
\newcommand{\qdd}{\Ddot{\bm{q}}}

\newcommand{\qs}{q}

\newcommand{\M}{\mathbf{M}}

\newcommand{\C}{\mathbf{C}}

\newcommand{\fe}{\bm{f}^{(\text{e})}}

\newcommand{\fI}{\bm{f}}
\newcommand{\fIs}{f} 

\newcommand{\Ktg}{\mathbf{K^t}}

\newcommand{\Ko}{\mathbf{K}^{(1)}}
\newcommand{\Ktw}{\mathbf{K}^{(2)}}
\newcommand{\Kth}{\mathbf{K}^{(3)}}

\newcommand{\ko}{\text{K}^{(1)}}
\newcommand{\ktw}{\text{K}^{(2)}}
\newcommand{\kth}{\text{K}^{(3)}}

\newcommand{\Koi}{\mathbf{K}^{(1),i}}

\newcommand{\param}{\bm{p}}
\newcommand{\paramn}{\hat{\bm{p}}}
\newcommand{\Np}{n_p}

\newcommand{\qr}{\bm{\eta}}
\newcommand{\qrd}{\mathbf{\bm{\dot{\eta}}}}
\newcommand{\qrdd}{\mathbf{\bm{\ddot{\eta}}}}

\newcommand{\qrs}{\eta}

\newcommand{\fIr}{\mathbf{\tilde{f}}}

\newcommand{\Mr}{\mathbf{\tilde{M}}}
\newcommand{\Cr}{\mathbf{\tilde{C}}}

\newcommand{\Kor}{\mathbf{\tilde{K}}^{(1)}}
\newcommand{\Ktwr}{\mathbf{\tilde{K}}^{(2)}}
\newcommand{\Kthr}{\mathbf{\tilde{K}}^{(3)}}

\newcommand{\Kori}{\mathbf{\tilde{K}}^{(1),i}}
\newcommand{\Ktwri}{\mathbf{\tilde{K}}^{(2),i}}
\newcommand{\Kthri}{\mathbf{\tilde{K}}^{(3),i}}

\newcommand{\kor}{\tilde{\text{K}}^{(1)}}
\newcommand{\ktwr}{\tilde{\text{K}}^{(2)}}
\newcommand{\kthr}{\tilde{\text{K}}^{(3)}}

\newcommand{\V}{\mathbf{V}}
\newcommand{\VT}{\mathbf{V}^T}

\newcommand{\Vi}{\mathbf{V}^i}
\newcommand{\VTi}{\mathbf{V}^{iT}}

\newcommand{\eigv}{\bm{\phi}}

\newcommand{\Pder}[2]{\frac{\partial#1}{\partial #2}}

\newcommand{\Pderm}[3]{\frac{\partial^2 #1}{\partial #3 \partial #2}}

\newcommand{\tensCon}[2]{\cdot_{{#1}{#2}}}

\newcommand{\pset}{\textbf{S}_{\textbf{p}}}
\newcommand{\setVMsRB}{\mathbf{S}_{\phi}}

\def\onedot{$\mathsurround0pt\ldotp$}
\def\cddot{
  \mathbin{\vcenter{\baselineskip.67ex
    \hbox{\onedot}\hbox{\onedot}}%
  }}%
\def\cdddot#1{
  \mathbin{\vcenter{\baselineskip.67ex
    \hbox{\onedot}\hbox{\onedot}\hbox{\onedot}%
  }}%
}

\newcommand{\rom}{\textit{ROM}}
\newcommand{\setone}[1]{\mathbf{S}_{#1}}
\newcommand{\settwo}[2]{\mathbf{S}_{#1}^{#2}}
\newcommand{\ntr}{N_t}
\newcommand{\nval}{N_v}

\newcommand{\sscr}[1]{\textsuperscript{\textcolor{blue}{#1}}}

\newcommand{\eps}{\epsilon}



\usepackage{babel}

\journal{Computer Methods in Applied Mechanics and Engineering}

\begin{document}

\begin{frontmatter}{}

\title{Non-Intrusive Data-Free Parametric Reduced Order Model for Geometrically Nonlinear Structures}

\author{Alexander Saccani\corref{cor1}}
\ead{asaccani@ethz.ch}

\author{Paolo Tiso}

\address{Institute for Mechanical Systems,\\ ETH Z\"urich, \\ Leonhardstrasse 21, 8092 Z\"urich, Switzerland}

\begin{abstract}
We present a fully non-intrusive parametric reduced-order modeling (PROM) framework for geometrically nonlinear structures subject to geometric variations. The method builds upon equation-driven Galerkin ROMs constructed from vibration modes and modal-derivative companion vectors, while nonlinear reduced tensors are identified from standard finite element outputs. A database of such ROMs is generated over a set of training samples, and all reduced operators-including the linear stiffness matrix, the quadratic and cubic nonlinear tensors, the Rayleigh damping parameters, and the reduction basis-are interpolated using Radial Basis Functions (RBFs). A global reduced basis is obtained through a two-level POD compression, combined with a MAC-guided reordering strategy to ensure parametric smoothness. The resulting PROM preserves the symmetry and polynomial structure of the reduced equations, enabling robust and efficient adaptation to new parameter values. Analytical parameter sensitivities follow directly from the interpolation model. The approach is demonstrated on a parametrically curved panel and a wing-box with geometric variations, showing excellent agreement with high-fidelity simulations and enabling substantial reductions in computational cost for parametric analyses.
\end{abstract}


\end{frontmatter}{}




\section{Introduction} \label{sec:Introduction}
Dynamic simulations of lightweight components play a pivotal role in the design of aerospace structures, allowing optimal design and  fatigue life prediction \cite{Gordon2011,Spottswood2008,Przekop2007,przekop2006nonlinear,mignolet2013review}.
At the foundation of these simulations lies the FE method that is used to derive a set of time-dependent ordinary differential equations whose time integration provides the prediction of structural displacements and stresses to the external loads.
The solved equations are nonlinear, as they are derived assuming a nonlinear strain-displacements relationship.
This relationship is either based on the \textit{Green-Lagrange} strain tensor or on the \textit{Von K\'arman} strain-displacement law \cite{crisfield1991}.
Moreover, the discretized geometries are usually complex, and they require fine meshing to correctly capture the stress gradients. This leads to high dimensional models, referred to as \textit{High Fidelity Models} (HFMs).
The combination of high dimensionality and nonlinearity results in large solution times that can hinder the design of aerospace components, especially when repeated analysis are required.\\ \\
Equation-based \textit{Reduced Order Models} (ROMs) have been introduced in the few last decades to cut down the computational costs associated to dynamic analysis of structures.
The most classic approaches for ROM construction rely on the \textit{Galerkin Projection} (GP) method, whereby the solution to the HFM equations is searched in a linear subspace spanned by few vectors stored in a \textit{Reduction Basis} (RB).
Then, the HFM equations are projected on the RB.
In equation-driven ROMs, the RB is assumed based on the physics of the problem.
For linear dynamical systems, a set of low frequency \textit{Vibration Modes} (VMs) constitutes an optimal RB, as proven by \textit{Modal Analysis} theory \cite{geradin2015mechanical}.
In fact, the typical external excitation usually excites only low-frequency modes.
A good RB for nonlinear systems can be constructed by complementing the RB of the underlying linear systems with additional vectors that capture the \textit{bending-stretching} coupling triggered by geometric nonlinearities in thin-walled structures.
Among these approaches,  the method of \textit{Static Modal Derivatives} (SMDs) \cite{idelsohn1985load} \cite{Idelsohn1985} and of \textit{Dual Modes} (DMs) \cite{Kim2013} have gained popularity and have been demonstrated to reliably approximate the response of geometrically nonlinear structures subjected to dynamic loading.
\\
After RB identification, the construction of the ROM is completed by computing the ROM operators associated to the projection of the HFM equations. While the reduced mass and stiffness matrices can be swiftly computed from their HFM counterpart through projection on the RB, the projection of the nonlinear forces is usually approximated with a cubic polynomial in the generalized reduced coordinates \cite{mignolet2013review}.
The coefficients of this cubic polynomial expansion are stored in third and fourth order tensors.\\
If the FE method is based on the \textit{Total Lagrangian Formulation} \cite{crisfield1991} and if the adopted material is linear elastic, this approximation becomes exact.
This could theoretically allow for the computation of the reduced tensors via direct projection of their HFM counterpart \cite{mignolet2013review}, following a so called \textit{direct-method}.
However, in most commercial FE programs, the FE implementation is based on the \textit{Updated Lagrangian formulation} \cite{crisfield1991} or on  \textit{Corotational formulations}, and the decomposition of internal elastic forces in their quadratic and cubic components is thus not accessible.
For this reason, reduced nonlinear tensors are identified using the FE program as a black box.
A set of equations in the unknown nonlinear tensors entries is formed by reading the nonlinear restoring forces or tangent stiffness matrices corresponding to a set of imposed displacements.
The method using the nonlinear force goes under the name of \textit{Enforced Displacement} (ED) method \cite{Muravyov2003} , while the one exploiting the tangent stiffness matrix is known as \textit{Enhanced Enforced Displacement} (EED) method \cite{Perez2014}.
Since tensors are identified only by using common \textit{input-output} quantities, both these methods are said to be \textit{non-intrusive} \cite{mignolet2013review}.
\\  \\
Unlike methods where the RB is constructed from data compression of HFM simulations, using for example, \textit{Proper Orthogonal Decomposition} (POD), equation-driven ROMs are relatively computationally cheap.
However, the computational cost for model construction can still be significant, and is usually dominated by the cost for reduced tensor identification. 
Moreover, in structural optimization \cite{Song2018,Hill2016,Park2023}, uncertainty quantification \cite{Marconi2021}, or model adaptation scenarios, this overhead has to be sustained multiple times.
A way to reduced the computational cost for ROM construction is to resort to \textit{Parameteric ROMs} (PROMs).
These techniques exploit an \textit{offline-online} decomposition.
In this setting, significant computational resources are used in the offline stage to identify a ROM that ideally works for any parameter realization in a prescribed parameter set. 
Then, ROM adaptation to any parameter realization is performed \textit{on the fly} in the online stage.
The initial overhead for PROM construction is usually larger than that for constructing a ROM for a prescribed parameter value. 
However, for multiple evaluations, the PROM strategy pays off, as ROM adaptation is extremely fast.\\ \\
Parametric data-driven ROMs techniques have reached a high level of maturity, while equation-driven ROMs lag behind.
A pioneering method for PROM construction for structural dynamics was put out in \cite{Amsallem2009} and extended in \cite{Amsallem2011}.
There, the authors devised a way to interpolate ROMs constructed with a RB of VMs while preserving important properties of the ROM operators, such as symmetry and positive definiteness.
The method, however, was applied only to linear dynamical systems.\\
In \cite{Park2023} the authors presented a parametric model for investigating the variability of the response of geometrically-nonlinear structures subjected to uncertainty in the geometry and boundary conditions.
The reduction strategy pursued in the paper is the Implicit Condensation and Expansion (ICE) \cite{McEwan2001}.
In \cite{Marconi2020} and \cite{Marconi2021} the authors presented a parametric ROM constructed with a RB of VMs and MDs to swiftly adapt the ROM to geometric imperfections.
This model was shown to capture the strong dependence of the response to the geometrical parameters.
However, it is intrusive in nature, as it requires the modification of the FE formulation at the element level.\\  \\
In this work we present a fully non-intrusive Parametric Reduced Order Model (PROM) for predicting the nonlinear dynamic response of structures undergoing geometric variations. Differently from previous approaches such as \cite{Marconi2020}  and \cite{Marconi2021}, the proposed method does not require access to element-level quantities and can therefore be implemented directly on top of commercial FE software (e.g., ABAQUS, ANSYS, NASTRAN). This non-intrusive character is essential for industrial applications, where the FE formulation is typically inaccessible.
The construction of the PROM proceeds in two stages. First, a database of equation-driven ROMs is generated at selected parameter samples using vibration modes and modal derivatives to capture nonlinear geometric coupling. Second, all reduced operators-including the reduction basis, the linear stiffness matrix, the quadratic and cubic nonlinear tensors, and the Rayleigh damping parameters-are interpolated across the parameter domain using Radial Basis Functions (RBFs). This interpolation strategy enables instantaneous adaptation of the ROM to any new parameter value and provides analytical parameter sensitivities through direct differentiation of the interpolants.
To the best of our knowledge, this is the first non-intrusive PROM for geometrically nonlinear structures that interpolates all reduced operators, including the nonlinear reduced tensors, while preserving their intrinsic symmetry and the polynomial structure of the reduced equations. Additionally, analytical gradients of the ROM in the parameters can be easily obtained by differentiation of the interpolation formula. These features make the proposed PROM particularly attractive for uncertainty quantification, structural optimization, and parametric model adaptation. \\ \\
We structure the presentation of our method as follows. In Section 2, we set out the theory underlying equation-based ROMs for geometrically nonlinear structures. in Section 3, we describe the construction of the PROM. In Section 4, we eventually apply the PROM to two examples, namely a curved panel and a wing box with geometrical shape variations.
\section{Preliminaries} 

\label{sec:preliminaries}

\subsection{Problem Statement}
\label{sec:HFMeqs}
The starting point for the construction of the PROM herein proposed are the equations of motion for FE models of thin-walled structures with geometric nonlinearities and linear elastic material, which read
\begin{equation}
    \M(\param) \qdd + \C(\param)\qd + \fI (\q,\param) = \fe(t,\param),    
    \label{eq:eomF}
\end{equation}
where $\q, \qd, \qdd \in \R^{n}$ are the vectors of nodal displacements,velocities and accelerations, $\M, \C \in \R^{n \times n}$ are the mass and viscous damping matrices, $\fI(\q), \fe \in \R^n$ are the vectors of internal and external forces respectively, and $\param \in \R^{\Np}$ is the vector of parameters that is used to account for variations in the geometry, material, and boundary conditions.\\
The elastic internal forces are derived by assuming a nonlinear strain-displacement relation based on the \textit{Green-Lagrange} strain tensor or on the approximate  \textit{Von K\'arm\'an} strains \cite{crisfield1991}.
This leads to a cubic formulation that writes:
\begin{equation}
    \fI(\q,\param) = \Ko(\param) \cdot \q + \Ktw (\param) \cddot (\q \otimes \q) + \Kth(\param) \cdddot\   (\q \otimes \q \otimes \q )
    \iff
    \fIs_{i} = \ko_{ij}\qs_{j} + \ktw_{ijk}\qs_{j}\qs_{k} + \kth_{ijkl} \qs_{j}\qs_{k}\qs_{l},
    \label{eq:tensDecF}
\end{equation}
where $\Ko \in \mathbb{R}^{n\times n}$ is the linear stiffness matrix, $\Ktw \in \mathbb{R}^{n\times n \times n}$ and $\Kth \in \R^{n \times n \times n \times n}$ are respectively the tensors associated to the quadratic and cubic forces, "$\otimes$" is the dyadic product operator, and '$\ \cdot\ $','$\ \cddot\ $','$\ \cdddot\ \ $' denote respectively, single, double and triple contraction. 
 The tensors $\Ktw(\param)$ and $\Kth(\param)$ are commonly referred to as \textit{nonlinear stiffness operators}, while $\Ko(\param)$ is the linear stiffness.\\
 These operators are parameter dependent if geometrical and material variations are considered.
 Notably, the linear stiffness matrix and the nonlinear tensors are sparse and symmetric with respect to all axis, as they can be derived through partial differentiation of a quartic potential function \cite{rutzmoserThesis}. 
 For instance  $\ktw_{ijk} = \ktw_{ikj} = \ktw_{jik} = \ktw_{jki} = \ktw_{kij} = \ktw_{kji}$.
 Moreover, we assume that the linear stiffness matrix $\Ko$ is positive definite.
 This is equivalent to assuming that the undeformed configuration is stable.
 \\
 In the following, we also adopt the \textit{Rayleigh Damping} model \cite{geradin2015mechanical}, where the damping matrix $\C$ is defined as a linear combination of the mass and stiffness matrices, $\M$ and $\Ko$, as
 \begin{equation}
     \C(\param) = \alpha(\param) \M(\param) + \beta(\param) \Ko(\param),
     \label{eq:rayDamp}
 \end{equation}
 with parameter-dependent $\alpha, \beta \in \R$.\\
 \\
 The objective of this work is the construction of a PROM for the investigation of the dynamic response of the HFM model, in Eq.\ref{eq:eomF}, for parametric variations defined over a prescribed parameter set.
 Specifically, we define this set as
 \begin{equation}
     \pset := \{ \param \in \R^{\Np} |\  p^{\text{min}}_i \leq p_i \leq p_i^{\text{max}}, \ i = 1,...,\Np \},
     \label{eq:pSet}
 \end{equation}
 where $p^{\text{min}}_i \in \R$ and $p_i^{\text{max}} \in \R$ are respectively the minimum and maximum prescribed values for $p_i$.\\
\subsection{Equation Driven Galerkin ROMs for Geometrically Nonlinear Structures}
\label{sec:ROMconstr}
In this subsection we provide the reader with a condensed outline on equation-driven Galerkin ROMs for systems with geometric nonlinearities, as this is the foundation for the construction of the PROM here presented.
\subsubsection{Galerkin-Projection}
The GP is arguably the most classical and established technique for the construction a ROM. 
It consists in restricting the solution space of the HFM to a linear subspace spanned by a RB, $\V \in \R^{n\times m}$ (with $m \ll n$), and in solving the projected HFM equations for the reduced coordinate vector, $\qr(t) \in \R^m$, :
\begin{align}
    \q(t) \approx \V \qr(t),\\
    \VT\M\V\qrdd + \VT\C\V\qrd + \VT\fI(\V\qr) = \VT\fe.
\label{eq:redEqsGP}
\end{align}
The ROMs constructed using GP preserve the \textit{Lagrangian Structure}  \cite{Farhat2015} associated to the HFM, a desirable feature that guarantees the ROM to inherit the same stability characteristics of the time integration scheme applied to the HFM equations \cite{Farhat2015}.
\subsubsection{Equation-Driven Reduction Basis}
\label{sec:basis}
The RB used in Eq. \eqref{eq:redEqsGP} dictates accuracy and efficiency of the ensuing ROM: the larger the RB, the higher the expected level of detail, yet the lower the computational gain.
In this work, the RB is identified following an equation-driven approach, where all the information required for RB construction is inferred from the physics modeled in the equations. 
In this way, the RB construction process is extremely efficient compared to  data-driven methods, such as POD, that rely on HFM simulations. \\
Instead, a set of low frequency VMs $\eigv_i \in \R^n$ are computed by solving the generalized eigenvalue problem defined as
\begin{equation}
    (-\M \omega_i^2 + \Ko) \eigv_i = \zero,\ \ \ \text{for}\ \ 
    i = 1,...,n_{\phi}\ ,
    \label{eq:eigLin}
\end{equation}
where $\omega_i$ is the $i$\textsuperscript{th} angular eigenfrequency.
Of these vectors, only a few are selected based on the \textit{Modal Participation Factor} (MPF) of the external load \cite{geradin2015mechanical} and the load frequency bandwidth, and are used to form the RB. 
Throughout the paper, we identify the selected modes with their mode number, which we collect in set $\setVMsRB \subset \N $.
\\
A RB of VMs is expected to capture the response of the structure as it vibrates about its equilibrium position in a linear regime. 
However, as the displacements vibration amplitudes increase to a point where geometric nonlinearities are activated, the ROM with RB of VMs only predicts usually an overly stiff response. 
To correct for this effect, additional vectors are added to the RB, serving as additional DOFs to capture the hallmarks of the nonlinear response.
For geometrically nonlinear structures, SMDs and DMs represent a good choice of companion vectors.\\ \\
The SMDs are defined from the pairs of VMs retained in the RB.
Specifically, the SMD $\bm{\theta}_{ij} \in \R^n$ stemming from VMs $\bm{\phi}_i$ and $\bm{\phi}_j$ is defined as
\begin{equation}
    \Ko\bm{\theta}_{ij} = - \left. \Pderm{\fI(\q)}{\eps_i}{\eps_j} \right|_{\eps_i,\eps_j = 0} ,
    \label{eq:SMDdef}
\end{equation}
where $\q = \bm{\phi}_i\eps_i+\bm{\phi}_j\eps_j$. In this way, SMDs are retrieved by solving all the linear systems of equations for all the different combinations of pairs of VMs in the RB.
It is worth mentioning that the r.h.s. in  Eq. \eqref{eq:SMDdef} can be approximated using finite differences \cite{Morteza}. 
Namely, if a central difference scheme is employed, then
\begin{equation}
    \left. \Pderm{\fI(\q)}{\eps_i}{\eps_j} \right|_{\eps_i,\eps_j = 0} 
    =  \left.\Pder{\Ktg(\q)}{\epsilon_i}\right|_{\eps_i,\eps_j = 0} 
    \bm{\phi}_j
    \approx \frac{\Ktg{(h\bm{\phi}_i)}-\Ktg{(-h\bm{\phi}_i)}}{2h}\bm{\phi}_j,
    \label{eq:SMDsRHS}
\end{equation}
where $\Ktg \in \R^{n\times n}$ is the tangent stiffness matrix, and $h$ is a user-defined scalar perturbation parameter.
As a result, the SMDs can be computed non-intrusively from commercial FE programs.
A major limitation in the use of SMDs as companion vectors to VMs arises from their quadratic increase in number with the number of VMs retained in the RB, leading to potentially inefficient ROMs.
In fact, if $n_{\phi}$ is the number of VMs in the RB, the total number of ensuing SMDs is $n_{\phi}/2 + n_{\phi}^2/2$.
Selection strategies for SMDs have been proposed in \cite{tiso2011optimal}, reducing in this way the total number of vectors in the RB.\\
Although derived from the different perspective \cite{Idelsohn1985, idelsohn1985load}  of quasi-static modal analysis, it can be easily seen that SMDs appear in the second order Taylor series expansion of the nonlinear static solution problem
\begin{equation}
        \fI(\q^*) - \sum_{i \in \setVMsRB} \epsilon_i \Ko \eigv_{i} = \zero.
        \label{eq:NlStSol}
\end{equation}
Interestingly, the same equations lie at the foundation of the DMs approach \cite{Hollkamp2005,Radu2003,Kim2013}.
In fact, the DMs are obtained by solving the equations in \eqref{eq:NlStSol} with a numerical scheme, for imposed loads in the directions of single VMs and linear combinations of pairs of VMs. 
Proper Orthogonal Decomposition is then used to identify the optimal linear subspace that spans the set of computed nonlinear solutions, retrieving the DMs \cite{Kim2013}.\\
Both DMs and SMDs are capable of capturing the nonlinear effects triggered by the nonlinear coupling between the VMs in the RB and the ones left out.
In fact, this coupling is exercised in the construction of the RB (Eq. \eqref{eq:NlStSol}) using loads that are orthogonal to all VMs not included in the RB. 
The contribution of the VMs left out from the RB, to the static solutions of Eq. \eqref{eq:NlStSol}, differs from zero only because of the nonlinearity.  
In this way, when using DMs or MDs as additional companion vectors to low frequency VMs, it is implicitly assumed that the modal dynamic coupling is closely related to the modal static coupling. 
The successful employment of SMDs in \cite{Marconi2020,Marconi2021,Morteza,Sombroek2018,Wu2016,Andersson2023,Jain2018} and DMs in \cite{Radu2003,Kim2013,mignolet2013review,Perez2014} confirms the validity of this assumption. \\ \\
In this way, the RB generally writes 
\begin{equation}
    \V = [\mathbf{\Phi},\mathbf{\Theta}] \in \R^{n \times m},
    \label{eq:rb}
\end{equation}
where $\mathbf{\Phi} \in \R ^{n \times n_{\phi}}$ and $\mathbf{\Theta} \in \R^{n\times n_{\theta}}$ are the matrices containing, respectively, the low frequency VMs, and the SMDs or DMs.
While the first part of the RB is mass orthogonal by construction, the second part is mass orthogonalized to guarantee good numerical conditioning.
For this purpose, the Gram-Schmidt algorithm \cite{Morteza} can be used. 
\subsection{Tensorial Formulation of ROM}
\label{sec:tensForm}
Time integration of the ROM in the form written in Eq. \eqref{eq:redEqsGP} usually leads to small, if any, computational gains compared to the HFM.
This is attributable to the expensive evaluations of the nonlinear internal  forces and Jacobians required for the reconstruction of their reduced counterpart. 
In fact, the required number of operations scales with the size of the underlying FE model.
Significant computational gains are achievable by leveraging an \textit{offline-online} decomposition, where the ROM operators are computed once and for all in the \textit{offline} model construction stage and used \textit{online} for efficient evaluations \cite{rutzmoserThesis,Rutzmoser2017}.
In this setting, the reduced mass and damping matrices are precomputed by projecting their HFM counterpart on the RB, while the reduced elastic forces are formulated as a cubic polynomial of the reduced coordinates:
\footnote{ "$\tensCon{i}{j}$" denotes tensor contraction of the $i$\textsuperscript{th} dimension of the first tensor with the $j$\textsuperscript{th} dimension of the second tensor (e.g. for $\mathbf{A} \in \R^{k\times k \times k}$ and $\mathbf{B} \in \R^{k\times k}$, $\mathbf{A}\tensCon{2}{1} \mathbf{B}$ is equivalent to $\text{A}_{ijk}\text{B}_{jl}$), and we have adopted Einstein's notation.}
\begin{equation}
    \fIr(\qr) \triangleq \VT \fI(\V \qr) =  \Kor \cdot \qr + \Ktwr \cddot (\qr \otimes \qr) + \Kthr \cdddot\   (\qr \otimes \qr \otimes \qr )
    \iff
   \tilde{\fIs}_{i} = \kor_{ij}\qrs_{j} + \ktwr_{ijk}\qrs_{j}\qrs_{k} + \kthr_{ijkl} \qrs_{j}\qrs_{k}\qrs_{l},
    \label{eq:redForcTens}
\end{equation}
Here, the reduced linear stiffness matrix, $\Kor \in \R^{m\times m}$, and the reduced nonlinear stiffness tensors, $\Ktwr \in \R^{m\times m \times m }$ and $\Kthr \in \R^{m \times m \times m \times m}$, are defined as
\begin{subequations}
\begin{align}
    \Kor &= \VT \Ko \V \hspace{5cm} 
     &\kor_{ij} = V_{li}\ko_{lk} V_{kj}\\ 
    \Ktwr &= \left( \VT \Ktw \cdot \V \right)\tensCon{2}{1} \V \hspace{3cm}
    \iff 
    &\ktwr_{ijk} = V_{li} \ktw_{lrp} V_{pk} V_{rj}\\ 
    \Kthr &= \left( \left( \VT \Kth \cdot \V \right) \tensCon{3}{1} \V \right) \tensCon{2}{1} \V
    &\kthr_{ijkl} = V_{pi} \kth_{pqrs} V_{qj} V_{rk} V_{sl}. 
\end{align}
\label{eq:tensProj}
\end{subequations}
Notice that the reduced order tensors are generally not sparse, as their full-order counterparts, but they are still symmetric. This feature
allows for improved efficiency in reduced forces evaluations and limited storage requirements \cite{rutzmoserThesis}.
Similarly, the reduced tangent stiffness matrix, which is used in the assembly of the Jacobian of the residual in time integration, can be obtained by differentiating Eq. \eqref{eq:redForcTens}.
This yields an expression that is a quadratic polynomial in the reduced coordinates, as reported in Appendix.\\
The computation of reduced tensors using Eq.\eqref{eq:tensProj} is possible only if the FE full order tensors are available.
This is usually not the case when commercial FE programs are used.
Consequently, the reduced order tensors are identified using non-intrusive techniques (a.k.a. indirect methods) that leverage standard outputs from commercial FE software.
The ED \cite{Muravyov2003} and EED \cite{Perez2014} are the most popular approaches for tensor identification. 
While in the ED method the tensors are identified by computing the internal forces to a set of imposed displacements, in the EED method tensor identification is carried out from evaluations of the tangent stiffness matrix, thus improving efficiency \cite{Perez2014}. \\
Non intrusive tensor identification is generally an expensive procedure, and usually dominates the cost for the construction of equation based ROMs, especially when the number of the vectors in the RB is large. 
In \cite{Saccani2024} it was shown that decent speed-ups can be achieved by combining EED with hyperreduction techniques.




\section{Non-Intrusive Parametric Reduced Order Model}
Before proceeding with the description of the method for PROM construction, we introduce the normalized parameter vector $\paramn$.
Normalization is convenient, since the entries in the parameter vector $\param$, defined in Eq. \eqref{eq:pSet}, may potentially differ by dimensions and order of magnitude.
This makes the comparison between different parameter configurations challenging.\\
For this reason, we define the normalized parameter vector $\paramn \in \R^{\Np}$ from $\param$ as 
\begin{equation}
     \hat{p}_i = \frac{p_i - p_{i}^{\text{min}}}{p_i^\text{max} - p_i^\text{min}}, \ \ \ \text{for}\ i=1,...,\Np,
     \label{eq:normVect}
 \end{equation}
 After normalization, the entries of the parameter vector can vary from zero to one, and the set of interest $\setone{\param}$ (defined in Eq.\eqref{eq:pSet}) becomes a unitary hypercube centered at the origin. 
 In the following, we denote this set by $\setone{\paramn}$.\\ \\
The equation driven PROM that we herein present leverages an offline-online decomposition, where in the offline stage we construct a database of ROMs that are subsequently interpolated online using RBFs. 
More formally, for a given realization of the parameter vector belonging to the parametric set of interest $\paramn \in \setone{\paramn}$, we define the corresponding ROM as the collection of six objects:
\begin{equation}
    \rom(\paramn) := \{\V(\paramn),\Mr(\paramn), \Cr(\paramn), \Kor(\paramn), \Ktwr(\paramn), \Kthr(\paramn)\},
\end{equation}
In this definition, $\V(\paramn) \in \R^{m}$ is the RB constructed at $\paramn$, and the reduced order model operators $\Mr(\paramn),$ $ \Cr(\paramn), \Kor(\paramn), \Ktwr(\paramn), \Kthr(\paramn)$ are defined as the projection of their HFM parametric-dependent counterpart on the parametric RB.
In this way, the reduced mass, stiffness, and damping matrices read
\begin{equation}
    \Mr(\paramn) = \VT(\paramn)\M(\paramn)\V(\paramn),\ \ \  
    \Kor(\paramn) = \VT(\paramn)\Ko(\paramn)\V(\paramn),\ \ \ 
    \Cr(\paramn) = \VT(\paramn)\C(\paramn)\V(\paramn),
    \label{eq:romOpDef}
\end{equation}
while the nonlinear tensors, $\Ktwr(\paramn),  \Kthr(\paramn)$, are defined through Eqs. \eqref{eq:tensProj}.\\ \\
In this work, we generate some training parameter samples and  collect them in a training set $\settwo{\paramn}{\text{train}} := \{\paramn^i\}_{i=1}^{\ntr}$, with $\paramn^i \in \setone{\paramn}$.
Then, for each of the parameter samples we construct a ROM that we store in a ROM database, defined as
\begin{equation}
\settwo{\text{ROM}}{\text{train}} = \{\rom (\paramn^i)\ |\ \paramn^i \in \settwo{\paramn}{\text{train}} \}_{i=1}^{\ntr}.
    \label{eq:romdatabase}
\end{equation}
This set of ROMs is used to train a RBF interpolation model that allows us to swiftly compute new ROMs, for not-sampled parameter values.
Specifically, we construct RBF interpolants for all the individual ROM operators (i.e. RB, reduced mass and stiffness matrices, reduced tensors of nonlinear forces, and reduced damping matrix).\\ \\
In the following, we split the exposition into three parts: first, we set forth the construction of the ROM database; second, we illustrate the interpolation scheme; third we present the properties of the proposed interpolation method.
\subsection{Reduced Order Models Database Construction}
\label{sec:databaseConstr}
\subsubsection{Smoothness Requirement for ROM}
Before proceeding with the construction of the ROMs database, we should guarantee that the way we define the single ROMs allows for interpolation, i.e. their variation with the parameters is continuous. 
Indeed, any attempt in interpolating a discontinuous function with continuous interpolants might potentially yield  rather poor results.
To ensure continuity of the ROM with $\param$, we need to guarantee that the HFM operators (mass matrix, damping matrix and nonlinear forces) and the RB are continuously varying in $\param$.
Under these conditions, the ROM obtained by GP is continuous by construction.\\
For FE models of structures with geometric nonlinearities and linear elastic material, the continuity of the HFM is ensured if variations in $\param \in \setone{\param}$ result in continuous variations in the mesh nodal coordinates, boundary conditions, or material properties.
Conversely, the RB varies continuously across different parameters realizations only for some careful choice and consistent ordering of its vectors.
Based on these observations, we construct ROMs  using a global RB that becomes parameter-dependent after mass orthogonalization. 
As shown in the following, this RB meets the continuity requirement.
\subsubsection{ROM Database Construction Algorithm}
\label{sec:databaseConstrAlg}
The construction of the ROM training database is presented in the following seven points, summarized in Algorithm \ref{alg:constructROMdatabase}.
\begin{enumerate}
    \item The first step towards ROM construction consists of sampling the admissible parameter space $\setone{\paramn}$, generating a user-defined number $N_t$ of training points that are collected in the training set $\settwo{\paramn}{\text{train}}$.
    The algorithm presented here for PROM generation can be virtually used with any sampling scheme, provided that good coverage of the parameter set in ensured.
    Both structured and scattered grids can be used, as the interpolation scheme that we adopt is based on RBFs, and can handle both cases. 
    However, in the authors' experience, scattered grids generated with quasi-random sampling schemes, such as \textit{Latin Hypercube Sampling} (LHS) \cite{helton2003latin}, provide the best compromise between low number of samples and accuracy of the interpolation model.
    \item As a second step, we construct a FE model for each point in $\settwo{\paramn}{\text{train}}$, and compute all the VMs within the frequency range of interest.
    Then, for each constructed model, we select the set of VMs that are excited at linear level by the load. 
    This selection is carried out based on the \textit{Modal Participation Factor} (MPF)\cite{geradin2015mechanical} and on the load frequency content, as detailed in \cite{tiso20214,Saccani2024}. 
    The selected VMs are collected in a matrix $\mathbf{\Phi}^{i} \in \R^{n \times n^{i}_{\phi}}$, where the index $i$ refers to the training sample number, and $n^{i}_{\phi}$ is the number of selected VMs. Notice that in principle, $n^{i}_{\phi}$ can be different for the different training samples. 
    \item The third step is devoted to the computation of RB vectors  needed to tackle nonlinearities, namely $\mathbf{\Theta}$ in Eq. \eqref{eq:rb}.
    Again, for each point $i$ in $\settwo{\paramn}{\text{train}}$, we construct a different matrix $\mathbf{\Theta}^i \in \R^{n \times n_{\theta}^i}$, where $n_{\theta}^i$ is the number of additional vectors. Both SMDs or DMs can be used, but it is important to mention that the proposed PROM construction method is not restricted to this choice. 
    In fact, any other approach that allows one to identify a good RB for capturing nonlinearities could be equivalently employed at this stage.
    \item In the fourth step, the global RB is constructed from the set of previously computed local RB, $\mathbf{\mathbf{\Phi}}^i$ and $\mathbf{\mathbf{\Theta}}^i$, for $i = 1,...,N_t$.
    Specifically, two matrices $  \mathbf{\Phi_G}$ and $ \mathbf{\Theta_G}$ are constructed by collecting together all the column vectors in $\mathbf{\Phi}^i$ and $\mathbf{\Theta}^i$, for all the training points:
    \begin{equation}
        \mathbf{\Phi_G} = [\mathbf{\Phi}^1,...,\mathbf{\Phi}^{N_t}], \ \  \mathbf{\Theta_G} = [\mathbf{\Theta}^1,...,\mathbf{\Theta}^{N_t}].
        \label{eq:snapMatr}
    \end{equation}
    Here, $\mathbf{\Phi_G} \in \R^{n \times \sum_i n_{\phi}^i}$, and $\mathbf{\Theta_G} \in \R^{n \times \sum_i n_{\theta}^i}$. \\
    Then, the column vectors in $\mathbf{\Phi_G}$ and $\mathbf{\Theta_G}$ are unit normalized, before separately applying \textit{Singular Value Decomposition} (SVD), and obtaining
    \begin{equation}
        \mathbf{\Phi_G} = \mathbf{L_\Phi}\mathbf{\Sigma_\Phi}\mathbf{R_\Phi}^T,\ \ \mathbf{\Theta_G} = \mathbf{L_\Theta} \mathbf{\Sigma_\Theta}\mathbf{R_\Theta}^T,
    \end{equation}
    where $\mathbf{L_\Phi} $ and $\mathbf{L_\Theta}$ are the left singular vectors, $\mathbf{R_\Phi}^T$ and $\mathbf{R_\Theta}^T$ are the right singular vectors, and $\mathbf{\Sigma_\Phi}$ and $\mathbf{\Sigma_\Theta}$ are the diagonal matrices of singular values.\\
    Next, we define the global RB as
    \begin{equation}
        \V = [\mathbf{\tilde{L}_\Phi},\mathbf{\tilde{L}_\Theta}],
        \label{eq:glBasis}
    \end{equation}
    where $\mathbf{\tilde{L}_\Phi} \in \R^{n \times m_{\phi}}$ and $\mathbf{\tilde{L}_\Theta} \in \R^{n \times m_\theta}$ are the snapshot matrices containing respectively $m_\phi$ and $m_\theta$ left singular vectors, stored respectively in $\mathbf{L_{\Phi}}$ and $\mathbf{L_{\Theta}}$, and corresponding to the largest singular values.
    A suitable choice of $m_\phi$ and $m_\theta$ is usually made using a relative truncation energy criterion, based on which the number of retained modes is chosen so that the energy of the projected snapshot matrix is retained up to a defined threshold.
    For a generic snapshots matrix $\mathbf{A} \in \R^{n \times N_s}$, with $N_s \leq n$, the fraction of the energy retained in its reconstruction with $m \leq N_s$ left singular vectors can be written as
    \begin{equation}
        e_{m} = \frac{\sum_{i=1}^{m} \sigma_{i}^2}{\sum_{i=1}^{N_s} \sigma_{i}^2},
    \end{equation}
    where $\sigma_i$ are the singular values.\\
    In this setting, two different relative energy threshold $e_{\phi}$ and $e_{\theta}$ can be defined, resulting potentially in to different number of vectors in $\mathbf{\tilde{L}_\Phi}$ and $\mathbf{\tilde{L}_\Theta}$.
    Selecting appropriate energy truncation thresholds is critical for accuracy of the ensuing ROM and is left to the analyst.
    However, based on the authors' experience, energy threshold values greater than $0.99$ generally yield accurate predictions.\\
    At this point, two remarks can be made.
    In principle, a global RB could be constructed by stacking together all the vectors in $\mathbf{\Phi_G}$ and $\mathbf{\Theta_G}$ in a single matrix. 
    However, the PROM derived using this RB would be extremely inefficient, as its size would be equal to the sum of the degrees of freedom of the local RBs over the entire set of training points.
    Moreover, this RB would carry a lot of unnecessary repeated information, since the variations in $\mathbf{\Phi}^i$ and $\mathbf{\Theta}^i$ are expected to be relatively small across different parameter realizations.
    Based on these observations, we propose to compress the global RB with POD, applied separately to the 'linear part' and to the 'nonlinear part' of the RB.
    This latter choice is made to allow more freedom in the RB compression process as to what pertains the accuracy of the resulting global RB in approximating the linear and nonlinear components of the local RB  (i.e. $\mathbf{\Phi}^i$ and $\mathbf{\Theta}^i$).
    In fact, based on the authors' experience, for the same global RB size, the best results are obtained using a higher energy threshold for the truncation of the left singular vectors of the linear part,  rather than for the nonlinear part, i.e. $e_{\phi} \geq e_{\theta}$. \\
    Additionally,  the POD compression process, and consequently the extracted compressed RB, are significantly impacted by the scaling of the vectors in the snapshot matrix in Eq. \ref{eq:snapMatr}.
    In contrast to snapshots matrices constructed from displacements of HFM simulations, where the relative scaling of the different snapshots carries a physical meaning, the scaling of the vectors in $\mathbf{\Phi_G}$ and $\mathbf{\Theta_G}$ lacks of physical interpretation. Indeed, both VMs and SMDs are defined up to a multiplicate constant. 
    For this reason, we  unit normalize the columns of the snapshots matrices prior to POD compression.
    By doing so, we aim for a compression that is optimal in a global sense, as POD weighs all the different vectors equally.
    \item In the fifth step, we mass orthogonalize the global RB in Eq. \eqref{eq:glBasis}, for all points in $\settwo{\paramn}{\text{train}}$. 
    This operation is needed to guarantee a well conditioned ROM.
    Mass orthogonalization is performed by computing the VMs of the different ROMs in the training set, obtained by projecting the parametric mass and stiffness FE matrices onto the global RB.
    Specifically, for each realization $\paramn^i \in \settwo{\paramn}{\text{train}}$, we solve the generalized eigenvalue problem
    \begin{equation}
        (-\VT \M(\paramn^i) \V \omega_j^2 + \VT \Ko(\paramn^i) \V) \tilde{\bm{\phi}}_j^i = \zero,\ \ \ \text{for}\ \ 
    j = 1,...,m,\ \  m = n_{\phi} + n_{\theta},
    \label{eq:eigLinRom}
    \end{equation}
    where $\tilde{\bm{\phi}}_j^i \in \R^{n}$, and $\omega_j$ are respectively the $j$th VM and angular frequency of the ROM in reduced coordinates.
    Using the VMs, the mass-orthogonal RBs are then retrieved by performing the change of variables
    \begin{equation}
        \V^i = \V \tilde{\mathbf{\Phi}}^i,
        \label{eq:rb_loc}
    \end{equation}
    with 
    \begin{equation}
       \tilde{\mathbf{\Phi}}^i = [\tilde{\bm{\phi}}_1^i,...,\tilde{\bm{\phi}}_m^i],
    \end{equation}
    where again index $i$ refers to $\paramn^i \in \settwo{\paramn}{\text{train}}$.\\
    After mass-orthogonalization, the global RB inherits local features, since, for each parameter $\paramn^i \in \settwo{\paramn}{\text{train}}$, it is transformed using the local parameter-dependent mass and stiffness matrices. 
    However, the range of $\mathbf{V}^i$, corresponding to different $\paramn^i$, remains the same, as the column vectors are just recombined to diagonalize the linear part of the different ROMs.
    \item The vectors in the different RBs $\V^i$ vary continuously across the parameter samples in the training set, if consistently ordered.
    In fact, if $\M(\paramn)$ and $\Ko(\paramn)$ vary continuously with $\paramn$, as a consequence, also the VMs of the ROM, obtained by solving the eigenvalue problem in Eq. \eqref{eq:eigLinRom}, vary continuously with $\paramn$ \cite{bruls2007}.
    However, the ordering of the column vectors in each RB plays a crucial role when it comes to guarantee continuity, column by column.
    Ordering vectors based on the corresponding ROM angular frequencies (see Eq. \eqref{eq:eigLinRom}) is unwarranted and can lead to discontinuities in the RBs across the parameter set.
    In fact, mode veering and mode crossing phenomena can potentially appear \cite{Stephen2009,Heirman2011,bruls2007}, affecting the ROM VMs. 
    When mode veering occurs, the eigenfrequencies corresponding to pairs of different VMs approach each other, as the parameters of the systems are varied, and then veer away without crossing.
    Close to the veering point, the shapes of the two involved VMs change smoothly, and the modes morph into one another \cite{Giannini2016,Heirman2011}.
    In contrast, in mode crossing the eigenfrequencies of two different VMs coalesce for a critical value of the parameters and then separate again.
    In this case, the VMs retain their shape before and after the crossing, without interacting \cite{Giannini2016}.
    At the crossing point, the space associated to the two crossing VMs becomes degenerate, and every combination of the two VMs is also a VM of the system \cite{geradin2015mechanical}.\\ \\
    Reordering algorithms have been employed in the literature to track branches of VMs, both in numerical and experimental procedures \cite{kim2007,bonisoli2013crossing,bruls2007,bonisoli2013crossing,Giannini2016}.
    These algorithms rely on the \textit{Modal Assurance Criterion} (MAC)  \cite{allemang2003modal}, a normalized index that expresses the similarity of different mode shapes.
    In this work, we use a mass weighted MAC number to perform consistent reordering of the vectors in the set of RBs, $\setone{\V} := \{\V^i\}_{i=1}^{N_t}$, with $\V^i$ defined in Eq. \ref{eq:rb_loc}.
    Specifically, consider two RBs, $\V^I$ and $\V^J$, corresponding to the pair of parameter realizations $\paramn^I,\paramn^J \in \settwo{\paramn}{\text{train}}$.
    For each pair of column vectors $\bm{v}_i^I$ in $\V^I$ and $\bm{v}_j^J$ in $\V^J$,
    we define the mass weighted MAC number as
    \begin{equation}
        MAC^{IJ}_{ij} = \frac{   (\bm{v}_i^I\ ^T  \M^{(I)} \bm{v}_j^J)^2  }{  \bm{v}_i^I\ ^T  \M^{(I)} \bm{v}_i^I \ \cdot \ \bm{v}_j^J\ ^T  \M^{(I)} \bm{v}_j^J        }, \  \  \ \text{for} \ i,j = 1,...,m,
        \label{eq:MAC}
    \end{equation}
    where $\M^{(I)} = \M(\paramn^I)$ is the  mass matrix computed for parameter $\paramn^I$ \footnote{We could have also used $\M^{(J)}$. For practical purposes, nothing changes.}.
    The MAC number can vary from zero to one, indicating low modal correlation for values near zero and high modal correlation for values close to one.
   The mass weighted MAC number is preferred to the standard MAC number, as it is used to reorder mass orthogonal vectors - see  the definition of local RBs in Eqs. \eqref{eq:rb_loc} and \eqref{eq:eigLinRom}. 
    Based on the authors' experience, this choice results in a more robust reordering scheme compared to when the standard MAC number is used, as the difference in MAC numbers between correlated and uncorrelated pairs of vectors increases.\\ \\
    Reordering is performed one RB at a time. A reference RB  $\V^{I}$ is chosen from set $\setone{\V}$. Then, its vectors are related to the vectors of the RB to be reordered, here referred to as $\V^J \in \setone{\V}$.
    Specifically, the MAC number matrix is computed from the vectors in the two RBs using Eq. \eqref{eq:MAC}.
    Then, vectors in $\V^{I}$ are associated to the vectors of $\V^J$ that maximize the MAC number .
    Formally,  we can define a reordering vector $\bm{b}^{IJ} \in \R^m$ as
    \begin{equation}
        b^{IJ}_i  = \underset{1 \leq j \leq m}{\text{arg max}}\ MAC^{IJ}_{ij},\ \  \ \text{for}\ i=1,...,m,
        \label{eq:vecb}
    \end{equation}
    and the reordered RB $\tilde{\V}^J$ is obtained by reordering the columns of $\V^J$ as
    \begin{equation}
            \tilde{\V}^J(:,i) = \V^J(:,b^{IJ}_i),
    \end{equation}
    where we adopt index notation.\\ \\
    In the present Algorithm, vectors in different RBs are ordered using different reference RBs.
    This choice is motivated by the fact that RB vectors might exhibit strong variations with the parameters and using a constant reference might result in poor tracking of VMs braches.
    In fact, if the reference point far is too far away from the reorder point in the parameter space, correlation between the vectors in the two RBs might be lost.
    Consequently, for each RB that we order, we  choose as  reference RB, the one corresponding to the sample point in the parameter space that is the closest to the reordering point.
    The  adopted distance metric is the $l$-2 norm in the normalized parameter space:
        \begin{equation}
        d(\paramn^i,\paramn^j) = \| \paramn^i - \paramn^j\|_2,
        \label{eq:distP}
    \end{equation}
    where $\paramn^i$ and $\paramn^j$ are two realizations of the normalized parameter vector.\\
    We then define two different sets: one containing reordered RBs, referred to as $\setone{\tilde{V}}$, and one containing the RBs that still need to be reordered, namely $\setone{V}$.
    The idea is to subsequently draw RBs from the set of not ordered RBs $\setone{\V}$, reordering them and add them to the set of reordered RBs $\setone{\tilde{\V}}$, until all RBs are reordered.
    At each iteration, the selection of the RB to order from $\setone{V}$, and of the reference RB from $\setone{\tilde{V}}$ is performed based  on the distances between all the parameter values corresponding to the RB in the two sets.
    These distances are computed using Eq. \eqref{eq:distP} and stored in matrix $\mathbf{D} \in \R^{|\setone{\V}| \times |\setone{\tilde{\V}}|}$ as
    \begin{equation}
        \text{D}_{ij} = d(\paramn^i,\paramn^j).
    \end{equation}
    Here, $\paramn^i$ and $\paramn^j$ are parameter realizations whose associated RBs are respectively in $\setone{\V}$ and $\setone{\tilde{\V}}$.
     We then compare all the values in the distance matrix and select the RBs in $\setone{\V}$ and $\setone{\tilde{\V}}$ that correspond to the minimum distance in the parameter space.\\
    This translates in reordering RB $\setone{\V}(i^*)$ using the reference RB $\setone{\tilde{\V}}(j^*)$, where
    \begin{equation}
        (i^*,j^*) = \underset{1\leq i \leq |\setone{\V}|, 1\leq j \leq |\setone{\tilde{\V}}| }{\text{arg min}} \text{D}_{ij}.
    \end{equation}
    Loosely speaking, we define a starting reference point in the parameter space, and progressively move away from it, while reordering the nearby RBs by referring them to the closest RB in the reordered set, as shown in  Fig. \ref{fig:reordAlg}.
\begin{figure}
    \centering
\includegraphics[width=0.9\linewidth]{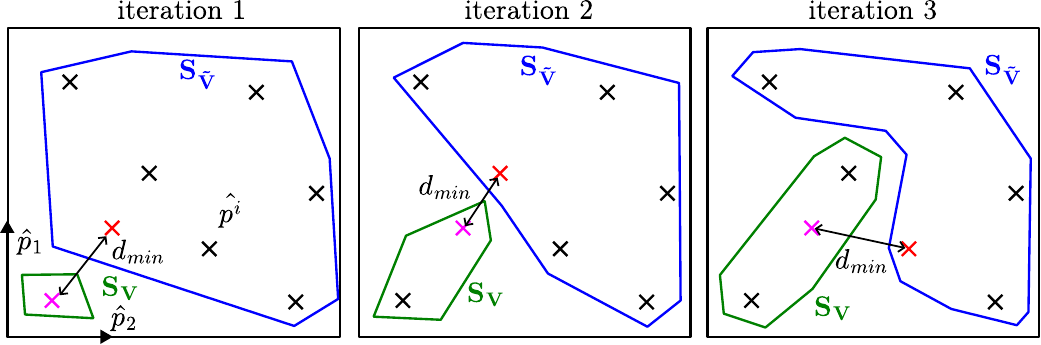}
    \caption{Schematic illustration of the reordering scheme used for RB reordering for a $2$D parameter space, for the first three iterations. At each iteration, we denote in magenta and red, respectively, the parameter points corresponding to the reference RB and to the RB that is reordered at the current iteration.}
    \label{fig:reordAlg}
\end{figure}
    The reader should note that an initial reference RB must be defined for the first reordering iteration. 
    This choice is arbitrary and is expected to affect the final ordering in the vectors of the different RBs but not reordering consistency.\\
    A limitation of RB reordering using the present approach is that two or more distinct vectors in the RB to be reordered, $\V^{J}$, may be associated to the same vector of the reference RB, $\V^I$.
    This situation arises when repeated values appear in the  reordering vector $\bm{b}^{IJ}$, defined in Eq. \eqref{eq:vecb}.
    As a result, the reordered RB becomes rank deficient and the ensuing ROM ill-defined.
    This issue highlights the inability of the MAC to distinguish different modes and is usually a symptom of an overly coarse parameter sampling, which causes too large variations in the shape of neighboring RBs vectors.
    An effective remedy is to refine the parameter sampling to reduce the distance between training points, thereby increasing the spatial correlation across RBs.
    \item The last step in the ROMs database construction consists of computing the reduced order operators from the full order operators using the orthogonalized, reordered RBs.
    In particular, reduced mass, stiffness, and damping matrices are obtained using Eq. \eqref{eq:romOpDef}, while the reduced nonlinear stiffness tensors are identified using either the ED or the EED method.
    \end{enumerate}
    Following the aforementioned procedure, we are able to construct a ROM database in the form presented in Eq. \eqref{eq:romdatabase}. 
    Moreover, the ROM operators are guaranteed to vary continuously across the parameter set of interest, allowing for interpolation.\\ \\
    We would like to stress that with the aforementioned procedure for ROM construction, we end up with a database of ROM whose linear part is diagonalized. 
    This proves to be a useful feature for ROM interpolation, as we will show in Section \ref{sec:ROMinterpProp}.
    Furthermore, if the ROM VMs (obtained using Eq. \eqref{eq:eigLinRom}) are mass-normalized, the ROM mass matrix $\Mr$ is the identity matrix, while the reduced stiffness matrix $\Kor$ is diagonal and contains the square of the ROM angular frequencies \cite{geradin2015mechanical}.
    Specifically, 
    \begin{equation}
    \tilde{\text{M}}_{jk}(\paramn^i) = \delta_{jk}, \ \ \tilde{\text{K}}_{jk}^{(1)}(\paramn^i) = \omega^{(i) 2}_j\delta_{jk}, \ \ \ \text{for} \ \ j,k= 1,...,m;
    \label{eq:redMassStiffDiag}
    \end{equation}
    with $\paramn^i \in \settwo{\paramn}{\text{train}}$, and $\delta_{ij}$ being the Kronecker delta.
    As a result, the Rayleigh damping reduced matrix is also diagonal and writes
    \begin{equation}
        \tilde{\text{C}}_{jk}(\paramn^i) = \alpha(\paramn^i) \delta_{jk} + \beta(\paramn^i) \tilde{\text{K}}_{jk}^{(1)}(\paramn^i)\delta_{jk} , \ \ \ \text{for} \ \ j,k= 1,...,m.
        \label{eq:redDampDiag}
    \end{equation}

    \subsection{Interpolation scheme for Parametric ROM}
    \subsubsection{RBF interpolant computation}
   The interpolation of the ROM database is carried out using RBF. 
   Let $\mathbf{Q}$ be any of the ROM operators $\V,\Cr,\Kor,\Ktwr,$ and $\Kthr$,
   \footnote{The mass matrix does not need to be interpolated, as it is constant and equal to the identity matrix across the parameter set, as shown in Section \ref{sec:databaseConstrAlg}}
   and $N_e$ be the total number of non-zero entries in $\mathbf{Q}$ that we collect here in vector $\bm{g} \in \R^{Ne}$.
   Any generic $j$th entry in $\mathbf{Q}$ is a real-valued function of the ROM parameter vector, $g_j:\paramn\in \R^{N_p} \mapsto \R$, and can be approximated using RBFs as
   \begin{equation}
       g_j(\paramn) \approx \sum_{i=1}^{N_t} \text{W}_{ji} \gamma(\| \paramn - \paramn^i\|_2),\ \ 
       \label{eq:interpRBF}
   \end{equation}
   where $\paramn^i \in \settwo{\paramn}{\text{train}}$ are the centers of the RBF interpolant, $\text{W}_{ji}$ the associated weights, and $\gamma: \R \mapsto \R$ is the RBF kernel function.
   This equation can also be written in matrix form as
   \begin{equation}
       \bm{g}(\paramn) \approx \mathbf{W} \cdot \bm{\gamma}(\paramn),
       \label{eq:RBFeval}
   \end{equation}
   where $\mathbf{W} \in \R^{N_e \times \ntr}$ is the matrix of weights and $\bm{\gamma} \in \R^{\ntr}$ with $\gamma_i = \gamma(\| \paramn - \paramn^i\|_2)$.
   The RBF weights are computed by evaluating the ROM operator vector $\bm{g}$ at the $\ntr$ training points in $\settwo{\paramn}{\text{train}}$ and solving the linear systems of equations for the matrix of weights $\mathbf{W}  \in \R^{N_e \times \ntr} $ that reads
   \begin{equation}
       \mathbf{G}^T = \mathbf{\Gamma}^T \cdot \mathbf{W}^T.
        \label{eq:weigthsRBF}
   \end{equation}
   Here, $\mathbf{G} \in \R^{N_e \times \ntr}$ and $\mathbf{\Gamma} \in \R^{\ntr \times \ntr}$ are defined as
   \begin{equation}
       \mathbf{G} = [\bm{g}(\paramn^1),...,\bm{g}(\paramn^{\ntr})], 
   \ \ 
   \mathbf{\Gamma} = 
        \begin{bmatrix}
            \gamma(\| \paramn^1 -  \paramn^1\|_2) && \gamma(\| \paramn^1 -  \paramn^2\|_2)  && \ldots && \gamma(\| \param^1 -  \paramn^{\ntr}\|_2)\\
            \gamma(\| \paramn^2 -  \paramn^1\|_2) && \gamma(\| \paramn^2 -  \paramn^2\|_2)  && \ldots && \gamma(\| \paramn^2 -  \paramn^{\ntr}\|_2)\\
            \vdots && \vdots && \ddots && \vdots \\
            \gamma(\| \paramn^{\ntr} -  \paramn^1\|_2) && \gamma(\| \paramn^{\ntr} -  \paramn^2\|_2)  && \ldots && \gamma(\| \paramn^{\ntr} -  \paramn^{\ntr}\|_2)\\ 
        \end{bmatrix}.
   \end{equation}
   The reader should note that the solution of these linear systems of equations is usually not computationally intensive, as it requires only one factorization of matrix $\mathbf{\Gamma}$, which is usually small.
   Moreover, interpolating ROM operators based on  Eq. \eqref{eq:interpRBF} is much faster than constructing a new ROM as described in Section  \ref{sec:ROMconstr}.
   In fact, the RBF model evaluation requires only $\ntr$ evaluations of the RBF kernel $\gamma(\paramn)$ and one matrix-vector multiplication.
    Then, after interpolation, the entries in the interpolated vector $\bm{g}$ are rearranged in the interpolated ROM operators $\Cr(\paramn), \Kor(\paramn), \Ktwr(\paramn), \Kthr(\paramn),$ and $\V(\paramn)$, retrieving $\rom(\paramn)$ in its canonical form.
    \\
    \subsubsection{Validation of RBF Kernel}
    In RBF interpolation, different kernels can be used, among which popular choices are Gaussian, multi-quadratic, and inverse multi-quadratic. 
    The kernel functions may feature a dependency on few hyperparameters that govern their shape.
    In the following, we denote the RBF parameter vector with $\bm{\epsilon}$.
    For example, the Gaussian kernel RBF, $\gamma_{g}(\bm{x},\epsilon) = e^{-(\epsilon\| \bm{x} - \bm{x}^i \|_2)^2}$, is a function of the shape parameter $\epsilon$ that governs its peakedness.\\
    The correct choice of RBF hyperparameters plays a crucial role on interpolation accuracy and \textit{ad-hoc} tuning can be made using a validation set.
    Specifically, different weights are computed using different RBF hyperparameters, resulting in different interpolants.
    Then, accuracy of the different interpolants, in approximating the function evaluations in the validation set, is compared, and the interpolation model associated with the lowest error is selected.\\ \\
    Following this scheme, we first construct different interpolants of the ROM operators, from the parameter samples in $\settwo{\paramn}{\text{train}}$, using different values of the hyperparameter vectors $\bm{\epsilon}$.
    This is done by repeatedly solving the linear systems in Eq. \eqref{eq:weigthsRBF}, obtaining $\bm{\epsilon}$-dependent weight matrices $\mathbf{W}(\bm{\epsilon})$.
    Then, we construct a validation parameter set $\settwo{\paramn}{\text{val}} := \{\paramn^i\}_{i=1}^{\nval}$, using the same sampling scheme adopted for the construction of the training set,
    and its associated ROM database $\settwo{\text{ROM}}{\text{val}} = \{\rom (\paramn^i)\ |\ \paramn^i \in \settwo{\paramn}{\text{val}} \}_{i=1}^{\nval}$. 
    This is done by following the procedure detailed in Section \ref{sec:databaseConstr}.
    Different approximations to the ROM operators are obtained as shown in  Eq. \eqref{eq:RBFeval} for all $\paramn^i \in \settwo{\paramn}{\text{val}}$, using weights $\mathbf{W}(\bm{\epsilon})$ corresponding to the different $\bm{\epsilon}$ values.
    An average relative error measure is defined as
     \begin{equation}
        e_{rel}(\bm{\epsilon}) = \sqrt{\sum_{i=1}^{\nval} \frac{\| \bm{g}(\paramn^i) - \mathbf{W}(\epsilon) \bm{\gamma}_{\bm{\epsilon}}(\paramn^i)  \|_2}
        {\| \bm{g}(\paramn^i)\|_2}}, \ \ \ \text{with} \ \paramn^i \in \settwo{\paramn}{\text{val}},
    \end{equation}
    where in vector $\bm{g}$ we collect the values of the exact entries of the ROM operators stored in  $\settwo{\text{ROM}}{\text{val}}$.\\
    Based on this error metric, we select the interpolation model corresponding to the value of $\bm{\epsilon}$ for which the error is the lowest.
    It is worth stressing that different RBF interpolation models are constructed separately for the entries of the different ROM operators, resulting in different optimal values of $\bm{\epsilon}$.
    This additional flexibility in the choice of the RBF hyperparameters allows us to improve the fit, as compared to the case where the same value of $\bm{\epsilon}$ is used for all operators.\\ 
    Moreover, we remind the reader that even if our interpolation model is validated here only with respect to RBF hyperparameters, in principle, one could use the validation set to optimally select the best RBF kernel function. 
    This could be achieved following the same procedure adopted for RBF hyperparameters validation, paying the cost of additional solutions of Eq. \eqref{eq:weigthsRBF}.
\begin{algorithm}[H]
\caption{ROM Database Construction}
\textbf{Input:} \texttt{constructFeModel}, \sscr{a} $\setone{\param}$, $\ntr$ \\ 
\textbf{Output}:  ROM database $\settwo{\text{ROM}}{\text{train}}$.
\label{alg:constructROMdatabase}
\begin{algorithmic}[1] 
    \Statex \textit{Generation of training set} 
    \State $\settwo{\param}{\text{train}} \gets \texttt{generateParamSamples}(\setone{\param},\ntr)$\sscr{b}
    \State $\settwo{\paramn}{\text{train}} \gets \texttt{normalizeParamSamples}(\settwo{\param}{\text{train}},\param_{\text{bounds}})$\sscr{c}
    \Statex \textit{Construction of the local RBs}
    \For{$i = 1,...,\ntr$}
    \State $\paramn^i \gets \settwo{\paramn}{\text{train}}(i)$ 
    \State $\texttt{FeModel}^i \gets \texttt{constructFEModel}(\paramn^i)$ \sscr{d}
    \State $\M^i \gets\texttt{FeModel.MassMatrix()}$,
    $\mathbf{K}^{(1),i} \gets\texttt{FeModel.LinStiffMatrix()}$
    \State $\mathbf{\Phi}^i \gets \texttt{computeAndSelectVMs}(\texttt{$\M^i,\mathbf{K}^{(1),i} $,$\fe(t)$})$ \sscr{e}
    \State $\mathbf{\Theta}^i \gets \texttt{computeBasisNL(\texttt{FeModel}, $\mathbf{\Phi}^i$)}$ \sscr{f}
    \EndFor
    \Statex \textit{Construction of the global RB}
    \State $  \mathbf{\Phi_G} \gets [\mathbf{\Phi}^1,...,\mathbf{\Phi}^{N_t}], \ \  \mathbf{\Theta_G} \gets [\mathbf{\Theta}^1,...,\mathbf{\Theta}^{N_t}]$
    \State  $ [\mathbf{L_\Phi}\mathbf{\Sigma_\Phi}\mathbf{R_\Phi}^T] \gets \texttt{SVD} (\mathbf{\Phi_G}),\ \ [\mathbf{L_\Theta}\mathbf{\Sigma_\Theta}\mathbf{R_\Theta}^T] \gets \texttt{SVD} (\mathbf{\Theta_G})$ \sscr{g}
    \State $\V \gets  [\mathbf{\tilde{L}_\Phi}(:,1:m_{\phi}),\mathbf{\tilde{L}_\Theta}(:,1:m_{\theta})]$ \sscr{h}
    \Statex \textit{Mass-orthogonalize RB}
    \For{$i=1,...,\ntr$}
    \State $\Mr^i \gets \VT \M^i \V,\ \ \Kori \gets \VT \Koi \V$
    \State $\tilde{\mathbf{\Phi}}^i \gets \texttt{VMs}(\Mr^i,\Kori)$ \sscr{i}
    \State $\V^i \gets \V \tilde{\mathbf{\Phi}}^i$
    \EndFor
    \Statex \textit{Reorder RBs}
    \State $[\V^1,...,\V^{\ntr}] \gets \texttt{reorderRBs}([\V^1,...,\V^{\ntr}],\settwo{\paramn}{\text{train}})$ \sscr{l}
    \ \ \ \ 
    \Statex \textit{Construct ROMs}
    \For{$i=1,...,\ntr$}
    \State  $\Mr^i \gets \VTi \M^i \Vi,\ \ \Kori \gets \VTi \Koi \Vi, \ \ \Cr^i \gets \VTi \C^i \Vi$
    \State $[\Ktwri, \Kthri] \gets \texttt{identifyNonlinearTensors}(\V^i,\texttt{FeModel}^i)$ \sscr{m}
    \State  $\settwo{\text{ROM}}{\text{train}}(i) \gets \{\V^i,\Mr^i,\Kori,\Cr^i,\Ktwri, \Kthri \}$
    \EndFor
\end{algorithmic}
\end{algorithm}
\vspace{-10pt}
\footnotesize
\noindent 
\sscr{a}\texttt{constructFeModel} is a function that constructs a FE model for each value of the parameter vector.\\
\sscr{b}\texttt{generateParamSamples} is a function samples the parameter set of interest. In this work we use Latin Hypercube Sampling.\\
\sscr{c} The training set is normalized as indicated in Eq. \eqref{eq:normVect}.\\
\sscr{d} The FE model corresponding to parameter $\paramn^i$ is constructed.\\
\sscr{e} \texttt{computeAndSelect} is a function that computes the VMs and selects $n_{\phi}^i$ of them based on the external excitation.\\
\sscr{f} \texttt{computeBasisNL} is a function that computes the part of the basis used to capture nonlinearities. Dual Modes and Modal Derivatives are used in this work. \\
\sscr{g} \texttt{SVD} function is used for SVD of the snapshot matrix. \\
\sscr{h} Here, the first $m_{\phi}$ and $m_{\theta}$ vectors are selected based on a user-defined truncation energy. \\
\sscr{i} Function \texttt{VMs} is used to compute the VMs of the single ROMs in the training set. \\
\sscr{l} \texttt{reorderBasis} function reorders the vectors in the RBs $\Vi$.  \\
\sscr{m} \texttt{identifyNonlinearTensors} function is used to identify nonlinear stiffness tensors, given the RB. The Enforced Displacement \cite{McEwan2001} or the Enhanced Enforced Displacement \cite{Perez2014} identification schemes could be employed.\\
\normalsize
\subsubsection{Properties of PROM interpolation}
\label{sec:ROMinterpProp}
In this Subsection we present the properties of the interpolation model herein proposed.
As reported in the literature, any good technique for constructing a ROM for structural dynamics should preserve the physics of the problem, to avoid incurring significant loss in accuracy, or worse, numerical issues pertaining to ROM time integration \cite{Farhat2015,Amsallem2011,Amsallem2009}.
Such properties are related to the \textit{Lagrangian Structure} inherent to the HFM, based on which:
\begin{enumerate}[label=(\roman*), itemsep=0pt, topsep=0pt]
        \item The mass matrix mass $\M$ is symmetric-positive definite, as the kinetic energy is always positive;
        \item The linear stiffness matrix $\Ko$ and the nonlinear stiffness tensors $\Ktw, \Kth$ are symmetric about all axes, as detailed in Section \ref{sec:HFMeqs}, since they are obtained by differentiation of an elastic potential function;
        \item The linear stiffness matrix $\Ko$ is positive definite, as the undeformed configuration is stable;
        \item The damping matrix $\C$ is positive definite, as the system is overall dissipative.
    \end{enumerate}
 A ROM is said to preserve the Lagrangian structure when the aforementioned properties are carried over from the HFM operators to their ROM counterparts.
 Based on this definition, Galerkin Projection ROMs are structure preserving, as shown in \cite{Farhat2015}.
 The interpolation technique herein presented is also structure preserving, with some caveats.
 In fact:
 \begin{enumerate}[label=(\roman*), itemsep=0pt, topsep=0pt]
        \item The reduced matrix mass matrix $\Mr $ is not interpolated, as it is constant and equal to the identity. Thus its symmetry and positive definiteness is guaranteed;
        \item The reduced linear stiffness matrix is diagonal $\Kor$ (see Eq. \eqref{eq:redMassStiffDiag}) and only its diagonal entries are interpolated,  preserving symmetry. 
        In contrast, its positive definiteness is not necessarily preserved. 
        Based on the authors' experience, if a positive definite RBF kernel is used and if the training grid is fine enough, the interpolated entries of $\Kor$ are positive.
        However, to avoid issues, we always check the sign of the entries in $\Kor$ after interpolation and before time integration.
        Since $\Kor$ is diagonal, if its entries are positive, then the matrix is positive definite.
        This allows us to confirm that this property is indeed preserved;
        \item The nonlinear stiffness tensors of the interpolated ROM, $\Kor $ and $\Kthr$, are symmetric by construction since, for each pair or triplet of symmetric entries, we interpolate their value and then store it in the tensors to respect symmetries.
        \item The reduced damping matrix $\Cr$ is obtained using Eq. \eqref{eq:redDampDiag}.
        Specifically, instead of interpolating it entry by entry, we interpolate the values of $\alpha(\paramn)$ and $\beta(\paramn)$.
        Eventually, we reconstruct $\Cr$ using those interpolated values and the interpolated $\Kor$.
        Again, an easy check on the positive definiteness of $\Cr$ can be made by looking at the signs of the interpolated $\alpha$ and $\beta$ damping parameters.
    \end{enumerate}
Once more, the positive definiteness of $\Kor$ and $\Cr$ is not enforced as in \cite{Amsallem2011,Amsallem2009}, but comes as a consequence of accurate interpolation.
In the investigated examples, these properties have always been satisfied; thus, there was no need to use a different, more involved, interpolation scheme.
However, we remind the reader that in critical cases where positive definiteness of $\Kor$ and $\Cr$ is difficult to attain, one could, in principle, use the interpolation strategy presented in \cite{Amsallem2009} for the linear part of the ROM, while keeping the ROM database construction procedure and the nonlinear tensor interpolation unmodified.


\section{Applications} \label{sec:applications}
In this section we present the construction of the parametric interpolated ROM for two different structures.
The first is a canonical curved panel, while the second is the  box of a NACA wing.
In both cases, parametric geometrical variations are considered.\\
Accuracy and computational performance of the parametric ROMs will be benchmarked against HFM simulations.
In this work, we used our \textit{Matlab} based in-house FE software \textit{YetAnotherFEcode} \cite{yafec}  as HFM and for constructing the parametric ROM.
Although our FE code allows for an intrusive construction of the ROM, we used it in a non-intrusive fashion to demonstrate that the presented methodology can in principle be applied using commercial FE programs, as demonstrated in \cite{Saccani2024},\cite{Morteza}.
\subsection{Curved Panel}
\subsubsection{Model description}
\begin{figure}
    \centering
    \includegraphics[width=1\linewidth]{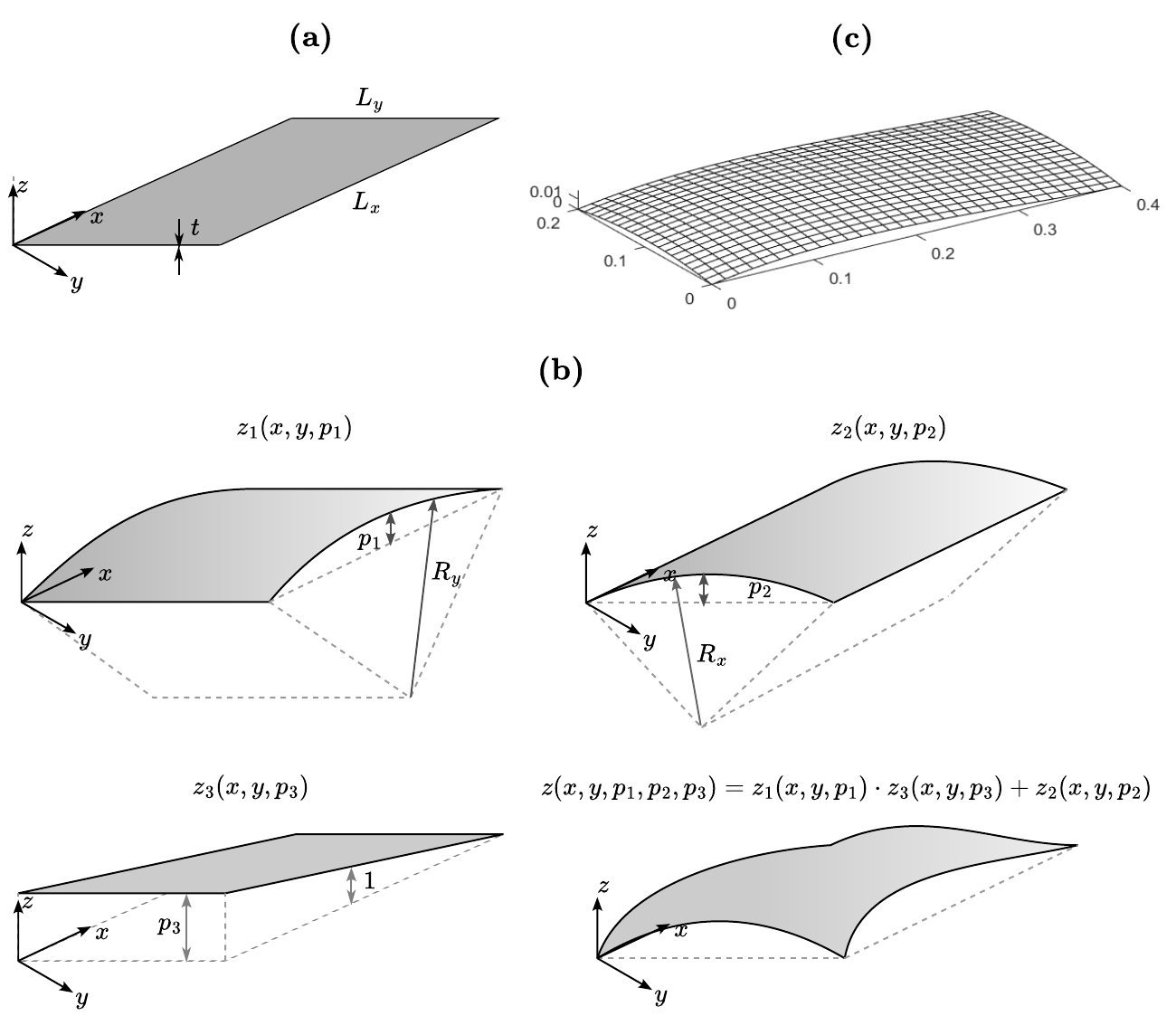}
    \caption{In (a) geometry of the flat panel. In (b) the parameters-dependent shape functions used to define the height of the plate. In (c), realization of the mesh for $p_1 = 1.3 \text{t}, \ p_2 = 0.89 \text{t}, \text{and}\  p_3 = 1.6$ (In the illustration, the out of plane $z$ coordinates are scaled by a factor of 10). }
    \label{fig:panelGeom}
\end{figure}

The rectangular curved panel under investigation is depicted in Fig.\ref{fig:panelGeom} a. It has dimensions $L_x = 40\  \text{mm}, L_y = 20\  \text{mm}$, thickness $\text{t} = 0.8 \ \text{mm}$, and is clamped at all four edges. 
The material is linear elastic with Young's modulus $E = 70\  \text{GPa}$, Poisson's ratio $\nu = 0.33$, and density $\rho = 2700 \ \text{Kg}/\text{m}^3$.
Three geometrical parameters, acting on the curvature of the panel, were considered.
Parametric variations can be used to model uncertainties introduced in the manufacturing process or, alternatively, to explore new designs.\\
In this test case, the final parameter-dependent geometry was obtained by superposing out-of-plane geometrical shapes to the flat rectangular panel, as shown in Fig..
Specifically, we constructed two out-of-plane geometrical shapes $z_1(x,y,p_1)$ and $z_2(x,y,p_2)$ by extruding two arch of circumferences passing through the edges of the plate along $y$ and $x$ axis, respectively.
Here, $p_1$ and $p_2$ are defined as the height of the midspan for each of the geometrical shapes.
The idea was then to add $z_1(x,y,p_1)$ to $z_2(x,y,p_2)$ to obtain the parametric out-of-plane shape.
In order to make the problem more challenging, we decided to introduce an additional parameter $p_3$ that skews the geometry towards one side of the panel, breaking symmetry.
For this purpose, we defined another shape function $z_3$ as $z_3(x,y,p_3) = p_3(1-\frac{2}{L_x})+\frac{2}{L_x}x$.
Here, parameter $p_3$ dictates the inclination of this plane that is forced to be equal to unity for $x = L_x/2$, for any value of $p_3$.
In this way, we defined the final out-of-plane shape of the panel as
\begin{equation}
    z(x,y,p_1,p_2,p_3) = z_1(x,y,p_1)z_3(x,y,p_3) + z_2(x,y,p_2).
\end{equation}
The parameter bounds considered here are shown in table \ref{tab:pboundsPanel}. 
As we show in the following, these parameter variations are enough to trigger significant variations in the panel response.
\\
\begin{table}[ht]
\centering
\caption{Parametric bounds for the curved plate}
\label{tab:pboundsPanel}
\begin{tabular}{|c|c|c|}
\hline
& $p_{min}$  & $p_{max}$ \\
\hline
 $p_1$& $0.250\text{t}$& $1.500\text{t}$ \\
\hline
 $p_2$& $0.620\text{t}$ & $1.875\text{t}$ \\
\hline
 $p_3$& $1$ & $2$ \\
\hline
\end{tabular}
\end{table}\\
Rayleigh damping model was assumed, where the mass and stiffness proportionality coefficients of the damping matrix, $\alpha$ and $\beta$,  were obtained by imposing a modal damping ratio of $0.01$ to the first and second linear VMs.
\subsubsection{Parametric ROM construction}
\begin{figure}
    \centering
    \includegraphics[width=0.4\linewidth]{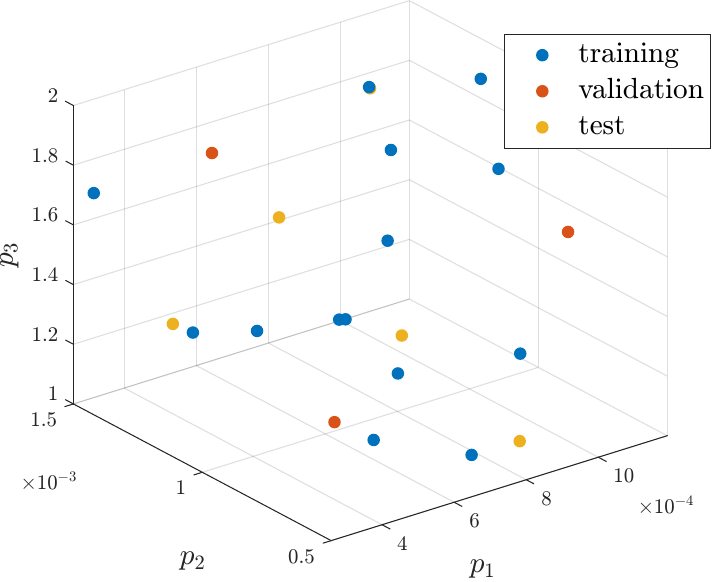}
    \caption{Training, validation and test samples of the curved panel.}
    \label{fig:samplingGrid}
\end{figure}

\begin{figure}
    \centering
    \includegraphics[width=1\linewidth]{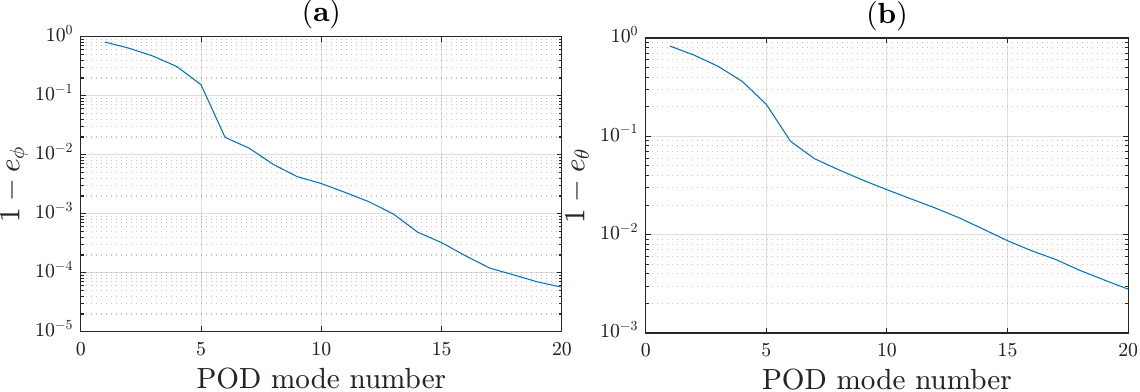}
    \caption{Energy lost in SVD data compression using different number of left singular vectors for $\mathbf{\Phi_{G}}$ (a) and $\mathbf{\Theta_{G}}$ (b).}
    \label{fig:PODenergy}
\end{figure}

The training set of parameters was constructed by drawing $\ntr = 14$ quasi-random samples using the Latin Hypercube Sampling scheme.
This sampling strategy maximizes the distance between points for given number of samples, ensuring a good coverage of the set of interest.
The generated samples are shown in Figure \ref{fig:samplingGrid}.
Then, the database of ROM was constructed following the procedure described in Section \ref{sec:databaseConstr}.
For each of the points in the training set, a FE model was constructed, using $600$ elements with quadratic shape functions, based on the \textit{Von-K\'arman} strain-displacement relationship \cite{crisfield1991}. 
All assemblies had the same topological equivalent mesh (i.e. same element connectivity) with $8505$ DOFs.\\
Then, for each FE model, the first 12 VMs were computed and selected based on the load Modal Participation Factor.
In this example, we considered the transient response of the panel to a uniform pressure load.  
Thus, we expected the response to be dominated at linear level by low frequency VMs.
For each assembly instance, we retained all the linear VMs in the frequency bandwidth [0,630 Hz] with non-zero MPF.
This led to the inclusion of VMs $1,2,3,8,10,11$ etc (after reordering had been applied).
The reader should note that the choice of the maximum frequency of interest was based on the analyst experience. 
\\
In this test case, we used the Static Modal Derivatives as complementary vectors used to capture nonlinearities. Their number grows almost quadratically with the number of retained VMs, hence it is essential to include only the most relevant one, to avoid undermining computational efficiency of the resulting ROM.
The Static Modal Derivatives were selected using the criteria proposed in \cite{Saccani2024}, which consists in ranking them based on the cross multiplication of the MPF of the different VMs. 
This is a variant of the more established selection criterion proposed in \cite{tiso2011optimal}, and has the advantage that Modal Derivatives selection comes with no additional computational cost (linear runs are not needed). 
For each of the assemblies, the selected SMDs were those coupling VMs $1-1,\ 1-2,\ 1-3,\ 2-2,\ 2-3,\ 3-3$. The SMDs were computed according to Eqs. \eqref{eq:SMDdef} and \eqref{eq:SMDsRHS}, where the tangent stiffness derivative was obtained with finite differences, using a  perturbation step $h = 10^{-8}$.
\\
The global RB was then constructed by unit normalizing all the VMs and Static Modal Derivatives, and stacking in them in two different snapshots matrices, $\mathbf{\Phi_G}$ and $\mathbf{\Theta_G}$, respectively.
Singular Value Decomposition was then applied separately to the two snapshot matrices, and 11 and 12 left singular vectors were selected from $\mathbf{\Phi_G}$ and $\mathbf{\Theta_G}$.
This corresponds to the relative truncation energy thresholds $e_{\phi} = 0.997$ and $e_{\theta} = 0.980$.
The energy lost in compression \footnote{$1-e$, where $e$ is the reconstruction energy} is plotted in Fig. for the two parts of the RB, $\mathbf{\Phi_G}$ and $\mathbf{\Theta_G}$, as a function of the number of retained POD modes.
As can be seen from Fig. \ref{fig:PODenergy} the graph has a kink in both two different plots in correspondence of the sixth POD mode.
This is because for each instance in the training set, six VMs and SMDs were chosen. 
The first six POD modes could be interpreted as an average of the first six modes throughout the training set.
Instead, the higher order POD modes are needed to capture parametric variations of these mode shapes across the training set.
\\
After the global basis was constructed,  $\ntr = 14$ different RB were obtained through mass normalization using the different mass matrices corresponding to the different FE instances in the training set.
Then using the different RBs, the reduced order mass, stiffness and damping matrices were obtained by projection of their parameter-dependent HFM counterpart, while the nonlinear tensors are identified with EED.
This required $299$ evaluations of the nonlinear tangent stiffness matrix.\\ \\
Interpolation of the ROM database was carried out in the normalized parameter space using RBF with inverse multiquadratic kernel function \cite{wright2003radial}, which writes
\begin{equation}
    \gamma(\delta,\epsilon) = \frac{1}{\sqrt{1+(\epsilon \delta)^2}}.
    \label{eq:invMulti}
\end{equation}
Here, $\delta = \| \paramn-\paramn^i\|_2$ is the $l_2$-norm distance and $\epsilon$ is a tuning parameter governing the peakedness of the distribution.
The best $\epsilon$ value was chosen with the aid of a validation set.
Specifically, $\nval = 3$ validation samples were used, and the validation ROM database was constructed following the same procedure for the construction of the training ROM database.
Different weights were computed for $50$ equally spaced $\epsilon$ values in a logarithmic scale, ranging from $10^{-2}$ to $10$.
Approximation accuracy, in reconstructing the ROM operators in the validation set, was compared for the different interpolants .
Based on this analysis, we selected the optimal $\epsilon$ values that we reported in Table \ref{tab:optValuesEps}. 
\begin{table}[ht]
\centering
\caption{Optimal $\epsilon$ values used in the interpolation of ROM operators of the curved panel.}
\label{tab:optValuesEps}
\begin{tabular}{|c|c|c|c|c|c|c|}
\hline
  & $\Kor$ & $\Ktwr$ & $\Kthr$ & $\V$ &$\alpha$ & $\beta$\\
\hline
 $\epsilon $& $0.012$& $0.534$ & $2.848$  & $0.201$  & $0.201$ & $0.231$  \\
\hline
\end{tabular}
\end{table}\\
\subsubsection{Parametric ROM interpolation performance}
\begin{figure}[h!]
    \centering
    \includegraphics[width=1\linewidth]{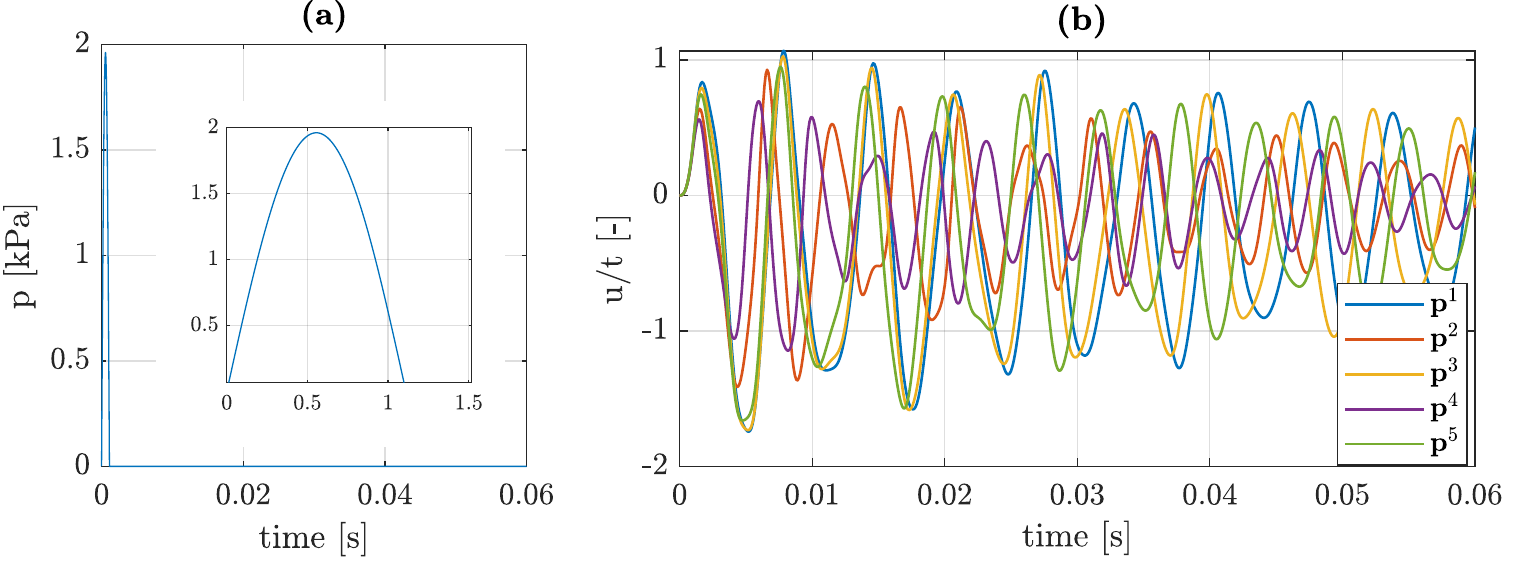}
    \caption{In (a) uniform pressure load time history, in (b) out-of-plane displacement time history of the midspan (normalized to the panel thickness) for the individual parameter realizations in the test set.}
    \label{fig:fullResp}
\end{figure}

\begin{figure}
    \centering
    \includegraphics[width=1\linewidth]{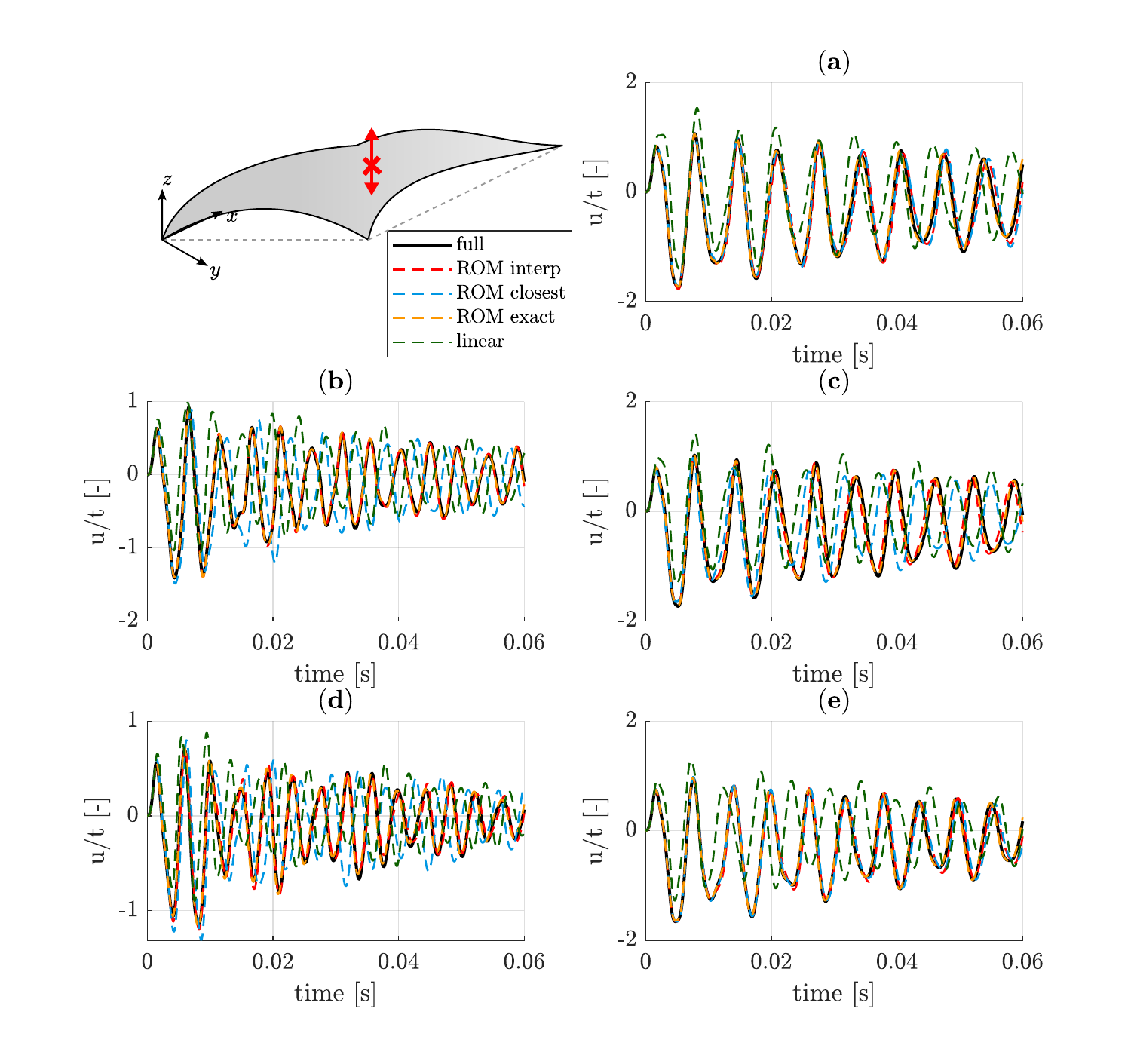}
    \caption{Out-of-plane displacement time history at panel midspan for different geometrical configurations in the test set.}
    \label{fig:outOfPlaneDisp}
\end{figure}

\begin{figure}
    \centering
    \includegraphics[width=1\linewidth]{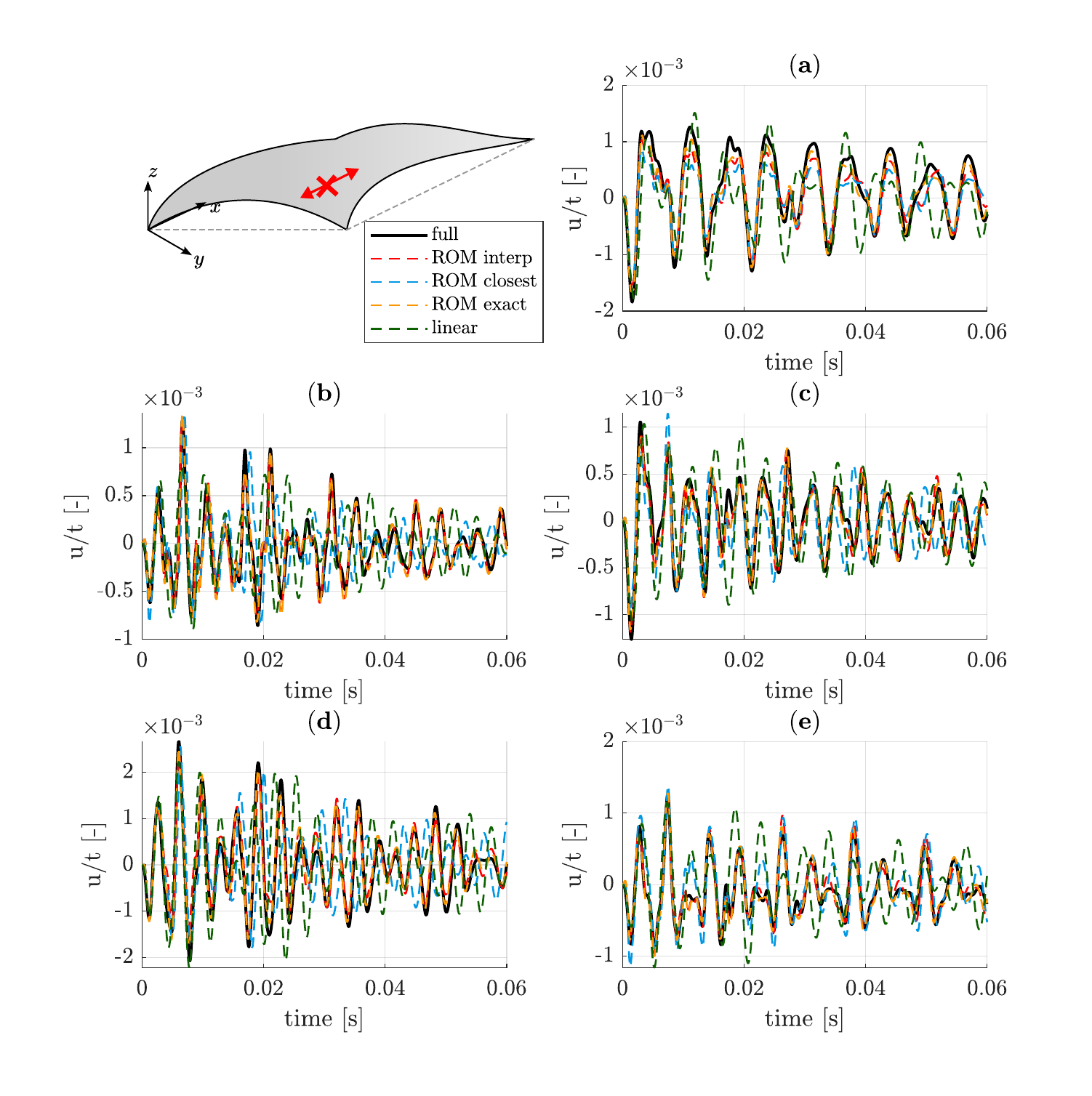}
    \caption{In plane displacement time history at panel midspan for different geometrical configurations in the test set. }
    \label{fig:inPlaneDisp}
\end{figure}

The accuracy of the interpolated ROM was evaluated using a test set containing five parameter samples that had not been used for the construction of the ROM.
These samples are shown in Figure \ref{fig:samplingGrid} and listed in Table \ref{tab:testSetVals} .
\begin{table}[ht]
\centering
\caption{Values of the parameter samples in the test set.}
\label{tab:testSetVals}
\begin{tabular}{|c|c|c|c|}
\hline
  $\param^i$ & $p_1$ & $p_2$ & $p_3$\\
 \hline
 $\param^1$& $1.09\text{t}$& $0.79\text{t}$ & $1.04$   \\
 \hline
 $\param^2$& $0.47\text{t}$& $1.61\text{t}$ & $1.32$   \\
 \hline
 $\param^3$& $0.82\text{t}$& $0.97\text{t}$ & $1.41$   \\
  \hline
 $\param^4$& $1.33\text{t}$& $1.84\text{t}$ & $1.77$   \\
   \hline
 $\param^5$& $0.54\text{t}$& $1.19\text{t}$ & $1.81$   \\
\hline
\end{tabular}
\end{table}
The test load was a short pulse modeled as a uniform pressure in space, acting in the positive $z$ direction, varying in time as
\begin{equation}
    p(t) = a\cdot\sin{(\frac{\pi}{T_{pulse}}t)}\cdot H(T_{pulse}-t),
    \label{eq:loadTime}
\end{equation}
where $H$ is the Heaviside, $T_{pulse} = 1/900 \ \text{s}.$ and $a = 1.962\  \text{kPa}$.
The time variation of the load is shown in Figure \ref{fig:fullResp} (a).\\
For each of the five test points, we run time integration over the time span $[0,0.06\ \text{s}]$, using five different models:
\begin{enumerate}[label=(\roman*), itemsep=0pt, topsep=0pt]
        \item the HFM, which is here the benchmark;
        \item the interpolated parametric ROM evaluated at the test points,
        \item the ROM in the training database corresponding to the training point that is the closest to the test point (we used as distance metric the $l_2$ norm in the normalized parameter space);
        \item the ROM obtained by projecting the global RB on the full order model evaluated at the test point, here referred to as recomputed ROM;
        \item the linearized model corresponding to the test point.
\end{enumerate}
In all cases, time integration was carried out using the \textit{Newmark-}$\beta$ integration scheme \cite{geradin2015mechanical}, with a constant time step $\text{d}t = 10^{-4} \ \text{s}.$
The panel responses for the out-of-plane degree of freedom (along $z$ axes) at the midspan obtained with the full model, for the different test parameter realizations, are shown in Figure \ref{fig:fullResp}.
A large variability of the transient response across the parameter set is observed, both in terms of frequency content and amplitudes of vibration.
This is also confirmed by the variation of the system first natural frequency, which ranges from a minimum value of $147.9\  \text{Hz}$ to a maximum value of $241.1\ \text{Hz}$.\\
A comparison of the results obtained from time integration of the different models is presented in Figures \ref{fig:outOfPlaneDisp} and \ref{fig:inPlaneDisp}. 
In Figure \ref{fig:outOfPlaneDisp} we plot the time history of the out-of-plane displacement for the node in the middle of the panel (same as shown in Figure \ref{fig:fullResp}), while in Figure \ref{fig:inPlaneDisp} the time history of the in-plane degree of freedom (along $x$ axes) of the node located at $x = 0.35L_x, y= 0.35L_y$.
For all points in the test set, it is possible to observe a good overlap of the time histories obtained with the HFM and the interpolated ROMs, both for out-of-plane and in-plane degrees of freedom.
The difference between the HFM with the linearized model confirms that the structure is vibrating in a nonlinear, as displacements in the order of the thickness also suggest.
For all test cases (Figure \ref{fig:outOfPlaneDisp} a - e), the characteristic period of the linear unforced vibrations appears to be smaller than the one of the nonlinear models, thus indicating a softening behavior of the nonlinear structure.
The comparison with the closest ROM in the training set is presented here to highlight the need for interpolation. 
Even if for some parameter realizations the interpolated ROM and the closest ROM show a similar behavior (Figure \ref{fig:outOfPlaneDisp} a and e), there are instances (Figure \ref{fig:outOfPlaneDisp} b,c, and d) where the difference is more marked, and using the closest ROM instead of the interpolated ROM leads to erroneous results.
Moreover, the reader should note that the recomputed ROM provides the best results, as the interpolated ROM is just its approximation. 
However, this comes at the price of reconstructing the ROM tensors, for each new parameter evaluation.\\
The average computational time to run time integration using the HFM and the interpolated ROM was, respectively, $848$ s and $23$ s, corresponding to a speed-up of $36.9$.
The reported computational times refer to a local machine equipped with $16 $ GB RAM and processor Intel Core i7-12700 H operating at $2.3$ GHz.
The modest computational gain is owing to the fact that the starting FE model is relatively small, especially compared to real-life industrial test cases, where the ROMs are clearly computationally advantageous.

\subsubsection{On ROM construction time}
Interpolating the database of ROMs to obtain a new adapted ROM is an extremely fast operation, compared to the computation of the tensors.
In fact, a new ROM can be constructed in less than a second once the interpolation weights are known, while the non-intrusive identification of tensors from the RB using the EED method requires on average around $1830$ s.
This second cost dominates the overall cost for model construction, as the SMDs and the VMs were computed in less than $10$ s, for each instance in the parameter set.
Tensor identification can be in principle parallelized, allowing for a major reduction in the construction time of the ROM database.
In this setting, the identification of tensors corresponding to different parameter instances can be dispatched to different workers.
The reader should note that the interpolation approach pays off only if multiple runs are needed, such as, for example, in uncertainty quantification or structural optimization settings.
In this second scenario, the interpolated ROM is also quite attractive because the gradients of the model with respect to the parameters, required for the search of the optimal parameter values, can be efficiently computed through analytical differentiation of the interpolation Equation \eqref{eq:RBFeval}.
If standard ROM techniques are employed, such gradients might be computed with Finite Differences or automatic differentiation techniques.
The first ones are computationally intensive, as they require multiple ROM constructions for values of the parameters perturbed about a reference point. 
The second ones are usually more computationally affordable and accurate, yet they require specific code and cannot be integrated in a non-intrusive ROM construction framework.

\subsection{NACA Wing Box}
\subsubsection{Model description}
\begin{figure}
    \centering
    \includegraphics[width=1\linewidth]{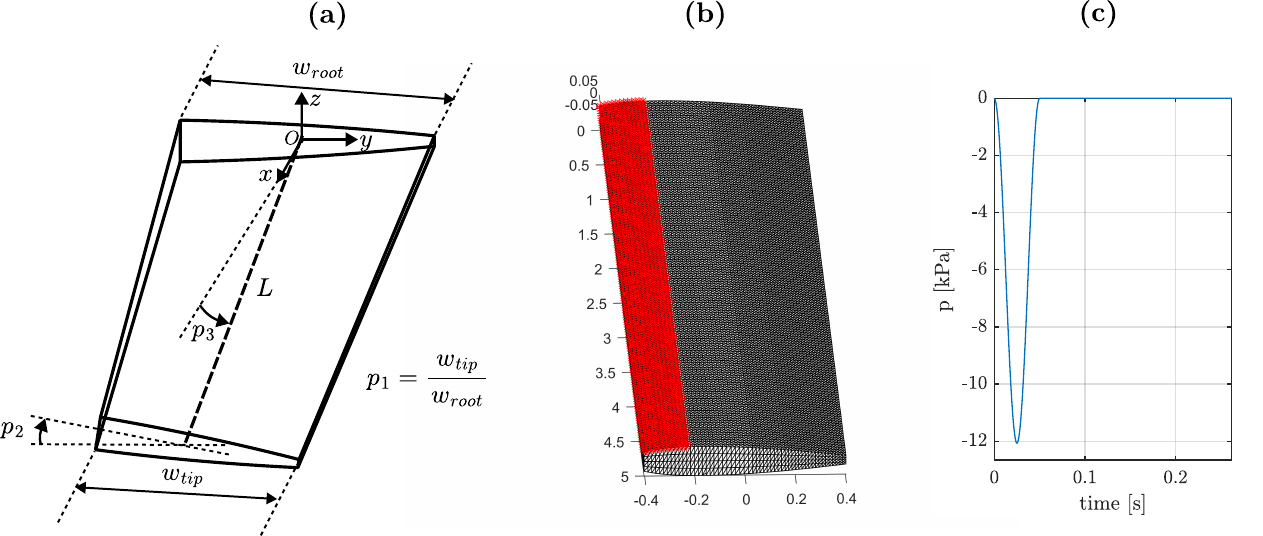}
    \caption{In (a) parametrized wing geometry. In (b) finite element model for zero tapering, twist and sweep angles. In (c) load time variation used for dynamic analysis. }
    \label{fig:wingModel}
\end{figure}

As a second test case we investigate the response of an aircraft wing with parametric geometry.  
The set of different geometrical configurations was defined starting from a nominal rectangular wing box with NACA 0012 profile. 
The wing consists of a skin panel stiffened with ribs along the longitudinal and lateral directions.
A model with the same geometry was considered in \cite{Jain2017, Jain2018}. 
Then the nominal wing geometry was morphed into different new geometric shapes by introducing three geometric parameters as shown in Figure \ref{fig:wingModel} a.\\
The first parameter, $p_1$, is the taper ratio and is defined as the ratio between the chord dimensions of the profiles at the tip and at the root of the wing. The profile sections across the span are rescaled by the taper ratio while preserving the proportions of the section.\\
The second parameter, $p_2$, is the angle of rotation of the tip wing cross section about the wing longitudinal axis passing from the origin of the centroid of the root section.
The third parameter, $p_3$, is the sweep angle which results in a inclination of the wing profiles towards the trailing edge of the wing.\\
The validity ranges for the parameters in the model are shown in Table \ref{tab:paramVarWing}.
The wing material is linear elastic with Young's modulus $E = 70\ \text{Gpa}$, Poisson's ratio $\nu = 0.3$, and density $\rho = 2700 \ \text{Kg}/\text{m}^3$. 
All sections in the model have a uniform thickness $\text{t} = 2\ \text{mm}.$\\ 
The FE model was constructed using $59040$ triangular shell elements with linear shape functions, featuring $3$ nodes per element and $6$ dofs per node.
The employed elements are based on the Von Karman strain-displacement kinematic relation.
Clamped boundary conditions constraining all the displacements and rotations at the root of the wing were applied.
This resulted in a model with $134195$ variables.
Meshes with same topology (same number of nodes, elements and connectivity) were used to model the different geometrical instances considered herein.\\
Objective of the study is the prediction of the displacement level of structures with different geometries to an external pulse pressure acting on the nodes highlighted in red in Figure \ref{fig:wingModel} b, along the negative $z$ direction. 
The load varies in time as in Eq. \eqref{eq:loadTime} with $T_{pulse} = 0.05\ \text{s}$ and its time history is plotted in Figure \ref{fig:wingModel} c.
We assume Rayleigh damping, with $\alpha$ and $\beta$ obtained by imposing a damping ratio of $1 \%$ to the first two structural VMs.
\begin{table}[ht]
\centering
\caption{Investigated parametric variations for the wing model}
\label{tab:paramVarWing}
\begin{tabular}{|c|c|c|}
\hline
& $p_{min}$  & $p_{max}$ \\
\hline
 $p_1$& $0.75$& $1$ \\
\hline
 $p_2$& $-3\  \text{deg}$ & $+3 \ \text{deg}$ \\
\hline
 $p_3$& $0$ & $+ 14.4\ \text{deg}$ \\
\hline
\end{tabular}
\end{table}\\
\subsubsection{Parametric ROM construction}
\begin{figure}
    \centering
    \includegraphics[width=1\linewidth]{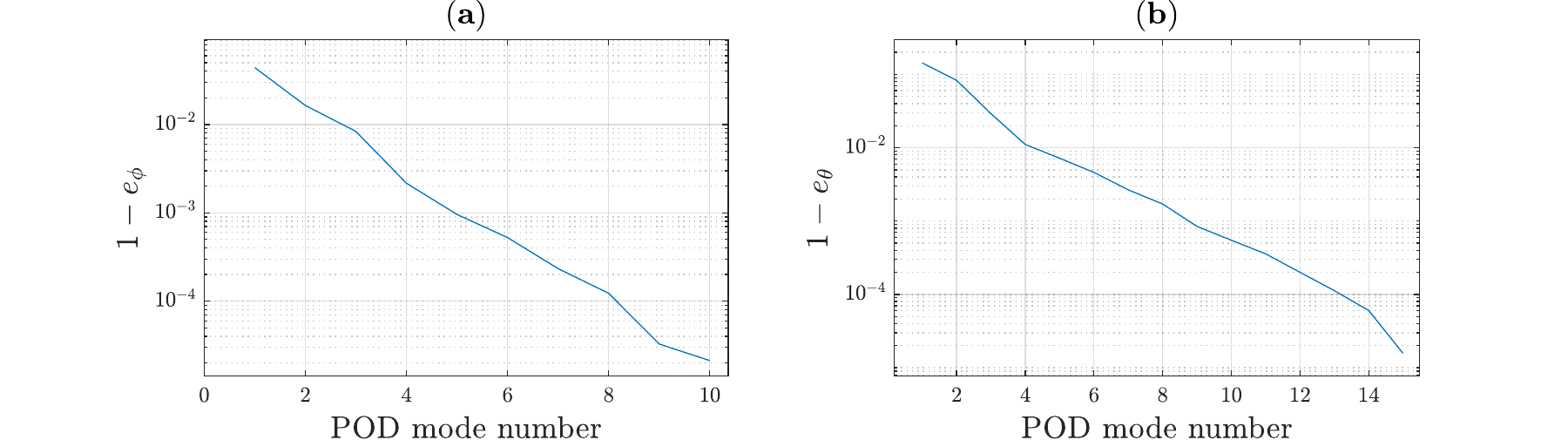}
    \caption{Energy lost in SVD data compression using different number of left singular vectors for $\mathbf{\Phi_G}$ (a) and $\mathbf{\Theta_G}$ (b), for the wing test case.   }
    \label{fig:podwing}
\end{figure}

Two different sets of $12$ and $4$ parameter samples were generated using Latin Hypercube Sampling and used as training and validation sets, respectively. 
Then, different FE models, corresponding to the different sampled parametric variations, were obtained by moving the nodes of the FE mesh corresponding to a wing with zero sweep and twist angles and null taper ratio.
For each FE model we computed the VMs and selected them based on the load Modal Participation Factor.
Based on this analysis we decided to construct a single-mode ROM, employing only the first VM, as the other modes excited by the load are much stiffer.
To prove the concept that the presented methodology can in principle be used with a generic RB for the nonlinear part (provided that it is representative of the nonlinear solution), we decided to use Dual Modes instead of SMDs.
Specifically, for each ROM we solved two nonlinear static problems for imposed loads in the direction of the first VMs, of the form
$\pm \Ko \phi_1 s_1$, where the $s_1$ is a scaling factor that is used to rescale the VM so that the maximum displacement of the mode at the tip of the wing is $3.5 \ \text{cm}$. 
This value was found to be sufficient to trigger the non-linearity while remaining in the neighborhood of the equilibrium point.\\
Then, VMs and Dual Modes corresponding to different models are unit normalized and stacked respectively in two different snapshot matrices $\mathbf{\Phi_G}$ and $\mathbf{\Theta_G}$, respectively.
Singular Valued Decomposition was applied separately to $\mathbf{\Phi_G}$ and $\mathbf{\Theta_G}$, and $5$ left singular vectors were extracted from each of the two matrices to form the RB for the parametric ROM.
This corresponds to a relative retained truncation energy threshold of $0.998$ for the VMs, and $0.99$ for the Dual Modes.
Energy decay of VMs and DMs with POD mode number are shown in Figure  \ref{fig:podwing}.
Then, this RB was mass-orthogonalized using the different mass matrices, corresponding to the different training samples, and the RB vectors reordering procedure was carried out.
Eventually, mass and linear stiffness matrices were computed by projection of their HFM counterparts, whereas nonlinear stiffness tensors were identified non-intrusively using the EED method. 
This last step required $65$ evaluations of the tangent stiffness matrix for each assembly in the training and validation sets.\\
Once the ROM database was constructed, RBF interpolation weights for the ROMs operators were computed for the inverse multiquadratic kernel in Eq. \eqref{eq:invMulti}.
The RBF shape parameter $\epsilon$ was tuned using the validation set. 
The $\epsilon$ values resulting from this procedure are reported in Table \eqref{tab:optValuesEpsWing}.
\begin{table}[ht]
\centering
\caption{Optimal $\epsilon$ values used in the interpolation of ROM operators of the NACA wing.}
\label{tab:optValuesEpsWing}
\begin{tabular}{|c|c|c|c|c|c|c|}
\hline
  & $\Kor$ & $\Ktwr$ & $\Kthr$ & $\V$ &$\alpha$ & $\beta$\\
\hline
 $\epsilon $& $0.0085$& $0.0077$ & $0.0085$  & $0.1385$  & $0.0077$ & $0.3854$  \\
\hline
\end{tabular}
\end{table}\\
\subsubsection{ROM interpolation performance}
\begin{figure}
    \centering
    \includegraphics[width=0.6\linewidth]{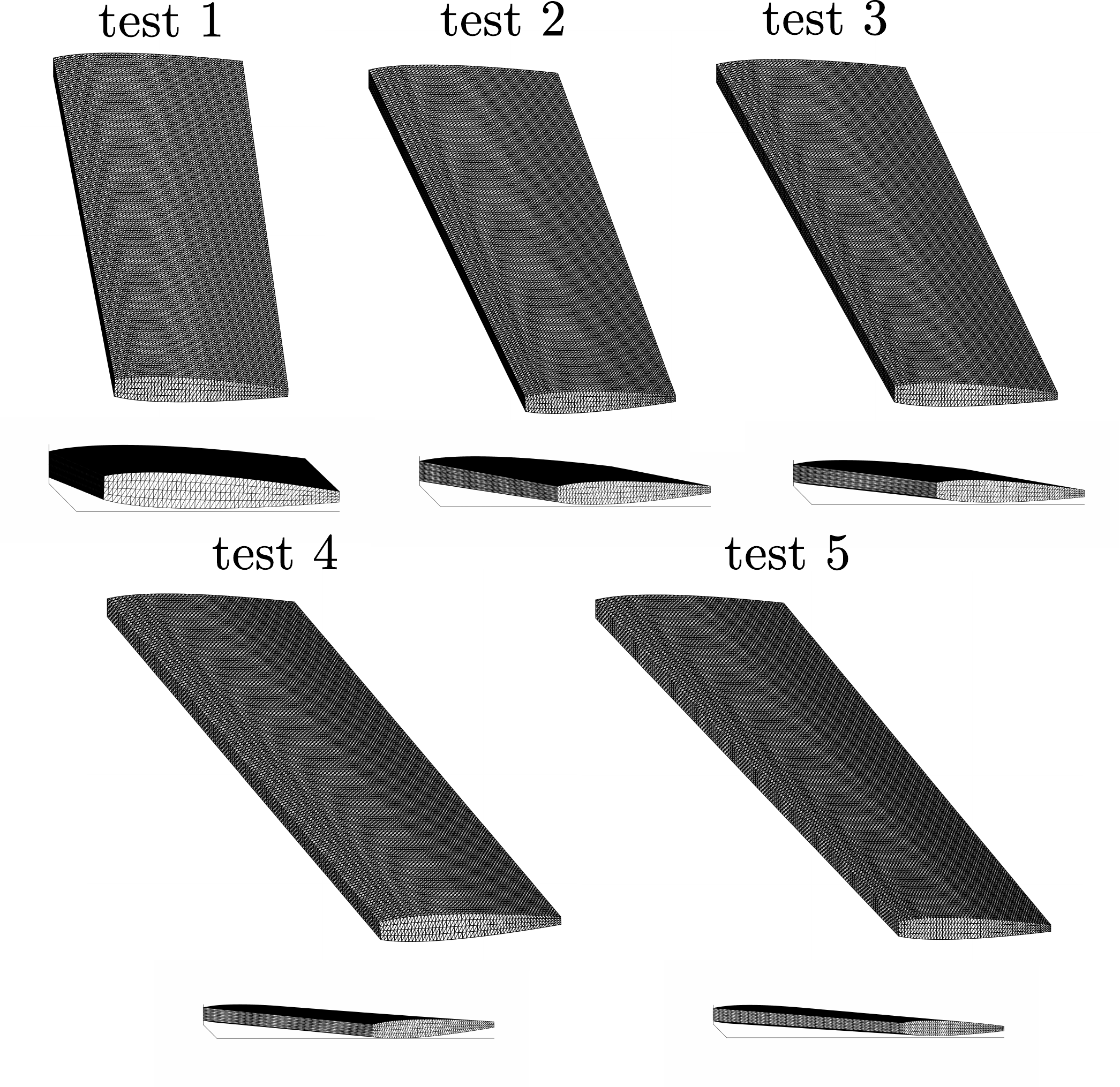}
    \caption{Geometrical configurations in the test set.} 
    \label{fig:testPoints}
\end{figure}

\begin{figure}[H]
    \centering
    \includegraphics[width=0.8\linewidth]{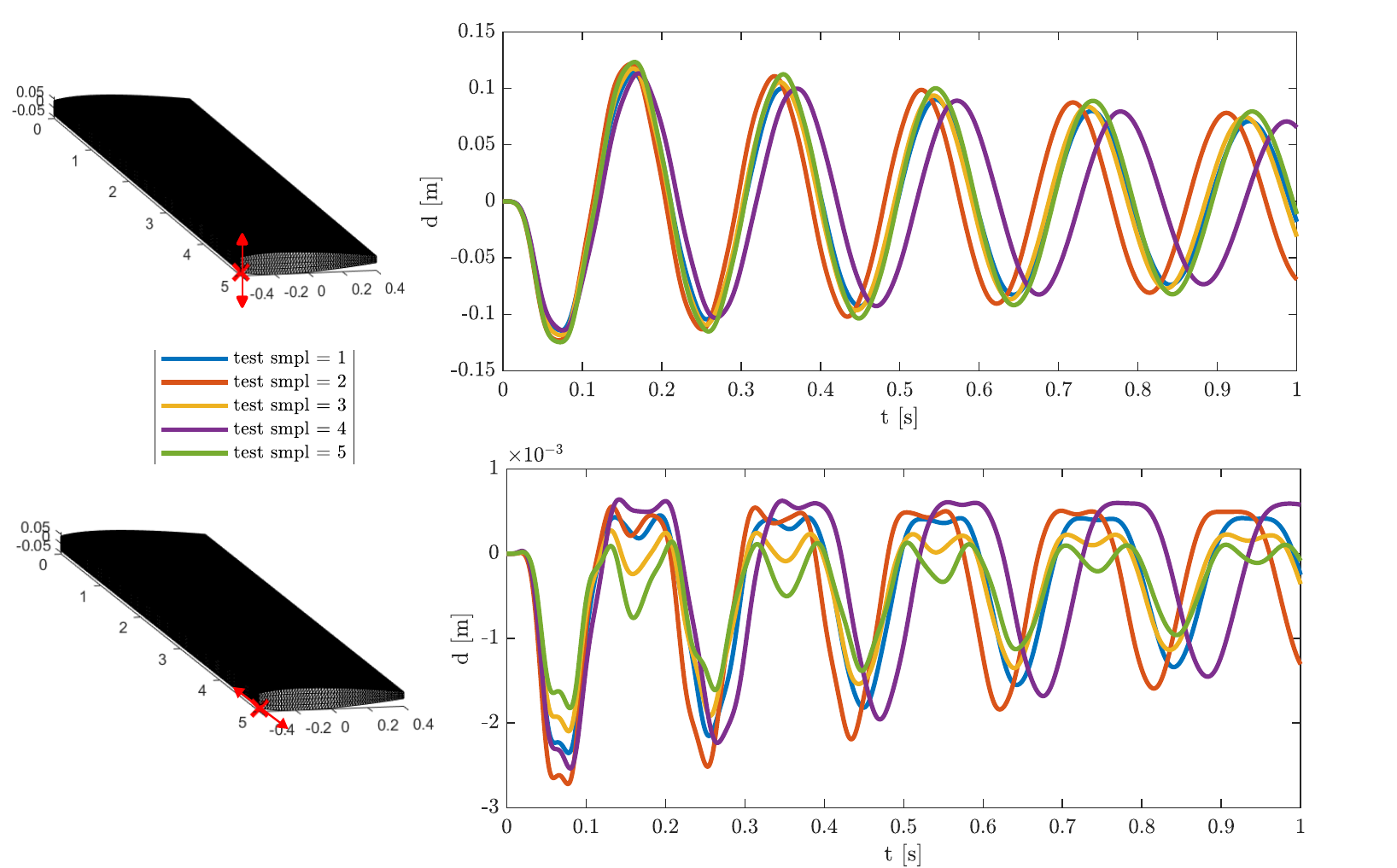}
    \caption{Wing out-of-plane (top) and in-plane (bottom) displacements of a node on the tip cross section for the different parameter realizations in the test set.} 
    \label{fig:fullWing}
\end{figure}

\begin{figure}
    \centering
    \includegraphics[width=0.8\linewidth]{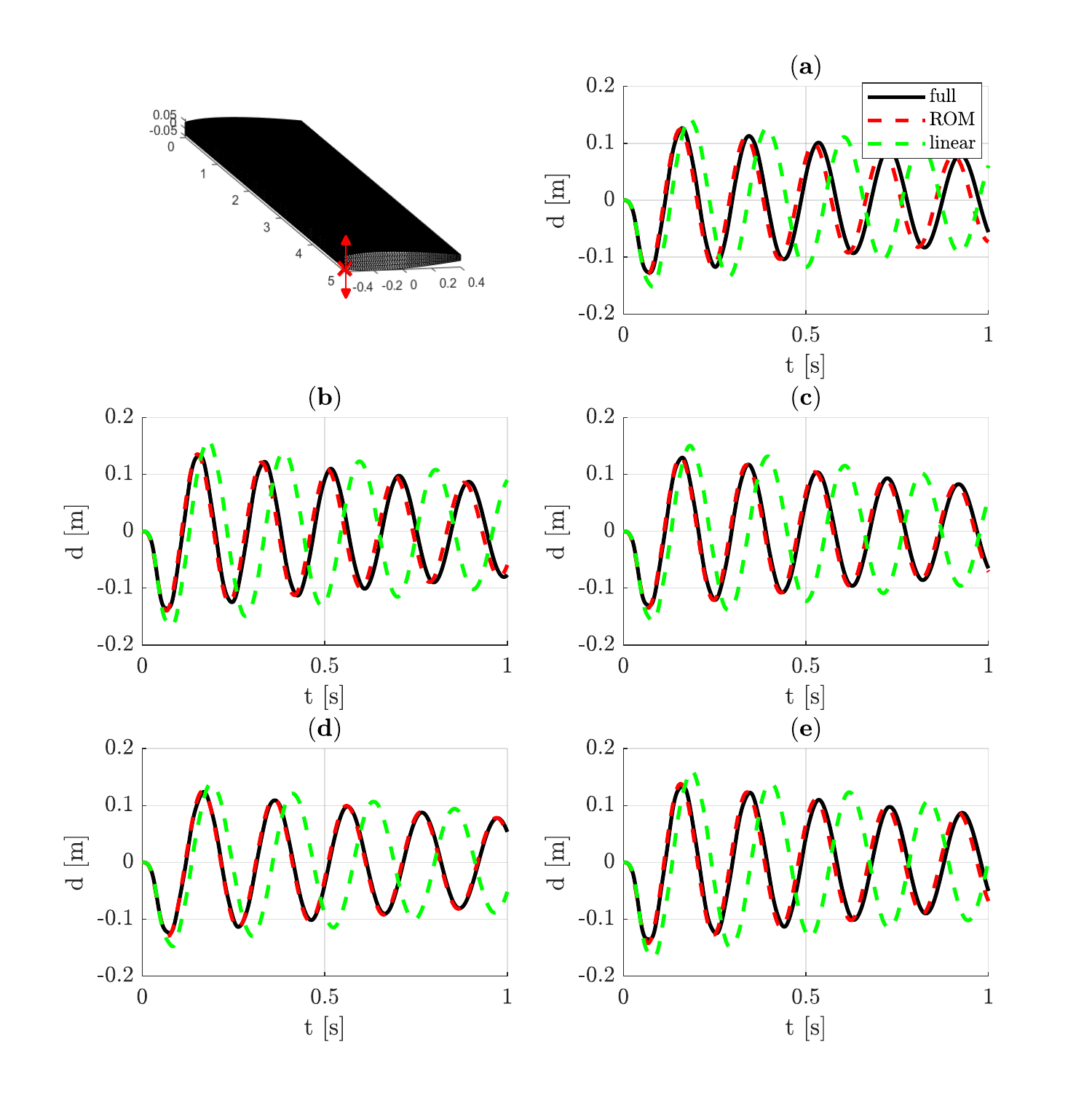}
    \caption{Out of plane response of the wing subjected to impulse loading at a node at the wing tip section. Comparison between HFM, interpolated ROM, and Linear model for the different parameter realizations in the test set, $\param^1$ (a), $\param^2$ (b), $\param^3$ (c), $\param^4$ (d), and $\param^5$ (e). } 
    \label{fig:OutOfPlane}
\end{figure}

\begin{figure}
    \centering
    \includegraphics[width=0.8\linewidth]{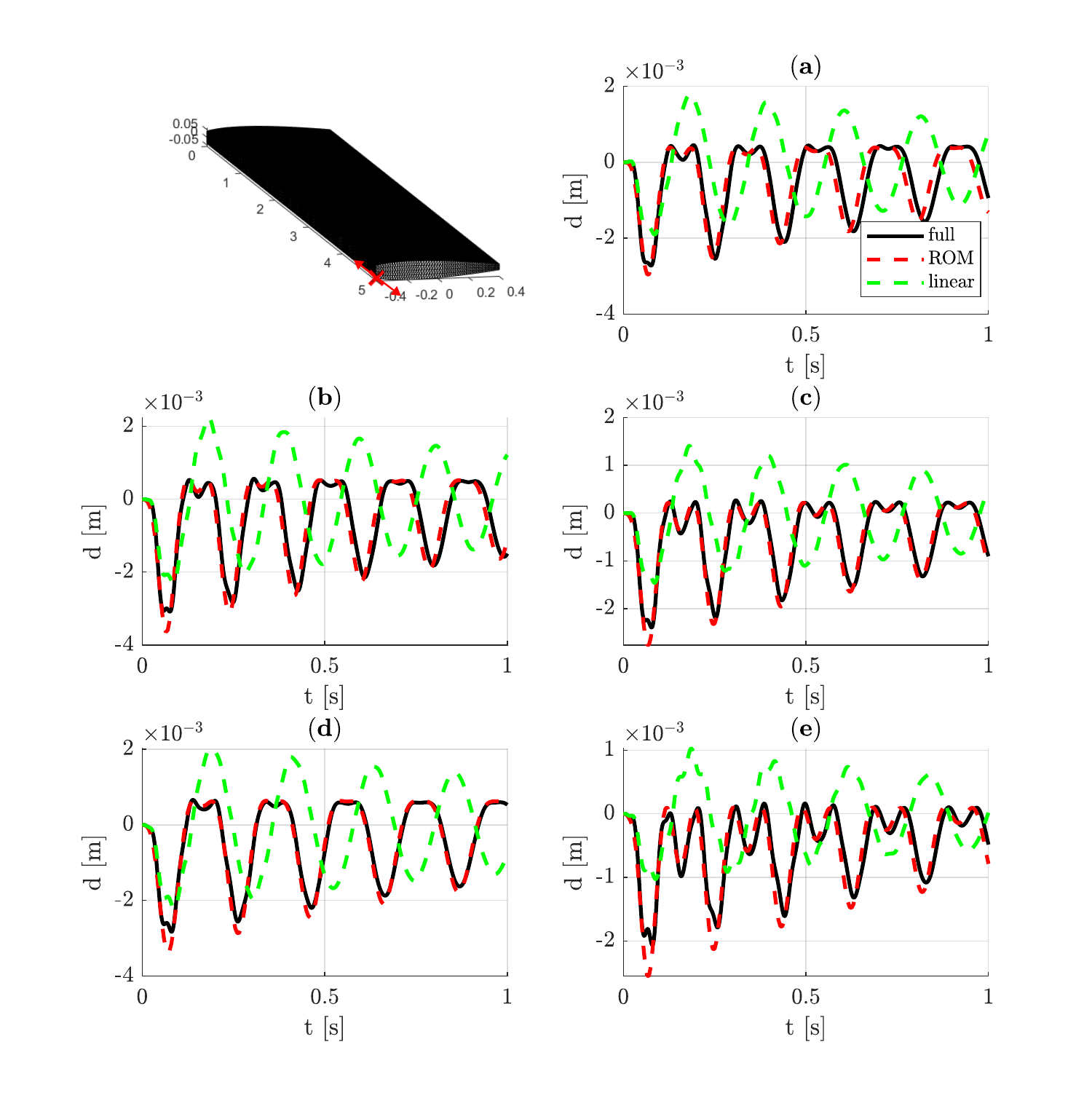}
    \caption{In plane response of the wing subjected to impulse loading at a node at the wing tip section. Comparison between HFM, interpolated ROM, and Linear model for the different parameter realizations in the test set, $\param^1$ (a), $\param^2$ (b), $\param^3$ (c), $\param^4$ (d), and $\param^5$ (e). } 
    \label{fig:InPlane}
\end{figure}

The parametric ROM performance is assessed by a direct comparison with HFM simulations for different random realizations of the parameter vector.
Specifically, a test set with five samples was constructed by sampling the feasible parameter space with LHS.
The values of the sampled parameters are reported in Table \eqref{tab:testSetValsWing}, while the corresponding geometrical configurations are shown in Figure \ref{fig:testPoints}.
\begin{table}[ht]
\centering
\caption{Values of the parameter samples in the test set for the wing.}
\label{tab:testSetValsWing}
\begin{tabular}{|c|c|c|c|}
\hline
  $\param^i$ & $p_1\ [\text{deg}]$ & $p_2\  [\text{deg}]$ & $p_3\ [-]$\\
 \hline
 $\param^1$& $0.92$& $-2.22$ & $0.61$   \\
 \hline
 $\param^2$& $0.80$& $1.73$ & $4.68$   \\
 \hline
 $\param^3$& $0.86$& $-1.33$ & $6.04$   \\
  \hline
 $\param^4$& $0.97$& $2.85$ & $11.15$   \\
   \hline
 $\param^5$& $0.81$& $-0.31$ & $11.80$   \\
\hline
\end{tabular}
\end{table}
For these configurations, time integration for the applied impulse load was performed using the Newmark scheme for:
\begin{enumerate}[label=(\roman*), itemsep=0pt, topsep=0pt]
        \item the HFM, which is here the benchmark;
        \item the interpolated parametric ROM evaluated at the test points;
        \item the linearized model corresponding to the test point.
\end{enumerate}
The displacements of a node on the left edge of the node tip section, obtained by time integrating the HFM for the parameter values in the test set are shown in Figure \ref{fig:fullWing}.
The reader should note that the variability in the parameters affects the characteristic frequency of vibration as well as the amplitude.
This variability is well captured by the interpolated ROM, as shown in Figures \ref{fig:OutOfPlane} and \ref{fig:InPlane}.
The difference between the linear response and the nonlinear prediction shows hardening behavior.\\

The HFMs were time integrated in parallel on the \textit{Euler} cluster of \textit{ETH Z\"urich}, using $5$ cores (one for different test assembly) and $4\ \text{Gb}$ RAM per core. 
No parallelization was used for the assembly of the internal forces and tangent stiffness matrices required within the time integration routine.
A constant time step of $1.2 \cdot 10^{-3} \text{s}$ was used. 
This resulted in an average total computation time of $4.789\  \text{h}$.
In contrast, the ROMs were time integrated using a larger constant time step $h = 5\cdot {10^{-3}}\ \text{s}$ and the total average solution time was $4\ \text{s}$.
This corresponds to a speed-up of the ROM with respect to the HFM of around $4310$.
The use of a larger time step to integrate the ROM compared to the HFM was possible as the high-frequency vibration effects, intrinsic to the HFMs, are not present in the ROMs.
\subsubsection{On ROM construction time}
The construction of the ROM was carried out in the Euler Cluster \textit{ETH Z\"urich}, using $16$ cores (one for different assembly in the training-validation sets) and $7\ \text{Gb}$ RAM per core. 
The most computationally intensive procedures in the construction of the database of ROMs needed for PROM construction were the following: 
\begin{enumerate}[label=(\roman*), itemsep=0pt, topsep=0pt]
        \item solving the nonlinear static problems needed in the Dual Modes computation, which required on average $62\  \text{s}$;
        \item identifying the nonlinear tensors using EED, which was completed on average in $120 \ \text{s}$;
\end{enumerate}
For this specific test case, the non-intrusive construction of the ROM is relatively cheap.
This owes to the fact that only few vectors are included in the RB.

\section{Conclusion} \label{sec:Conclusion}
In this work we presented a technique for constructing a PROM for structures with geometric nonlinearities, following an equation-based approach for the construction of the RB.
The PROM can handle parametric variations in the geometry, and potentially in the material and boundary conditions.
Thus, it can be used for model updating, uncertainty quantification, and structural optimization.\\
The PROM is constructed in two steps.
A database of ROMs is constructed for random parameters samples drawn from the parameter set of interest. 
Then, a RBF interpolant is fitted to the entries of the ROM operators to obtain the PROM.
he numerical examples demonstrate that the PROM reduces the computational effort for model adaptation to virtually zero, while still capturing substantial variations in the dynamic response across the parameter domain. Furthermore, analytical parameter sensitivities follow directly from differentiation of the PROM, a capability of particular relevance for gradient-based optimization.

Two limitations should be acknowledged. The first concerns the selection of an adequate number of parameter samples: too few lead to poor accuracy, while too many substantially increase the offline construction cost. Greedy sampling strategies \cite{Dubois} may provide an efficient remedy. The second limitation is inherent to the current geometric parametrization, which is restricted to topology-preserving mesh morphing. For example, in the context of stiffened panels, the method can account for variations in skin or stiffener thickness, but not for changes in the stiffener layout.


\bibliography{RRC2023_bib}

@misc{Amsallem2009,
   author = {David Amsallem and Julien Cortial and Kevin Carlberg and Charbel Farhat},
   journal = {INTERNATIONAL JOURNAL FOR NUMERICAL METHODS IN ENGINEERING Int. J. Numer. Meth. Engng},
   pages = {1-16},
   title = {A method for interpolating on manifolds structural dynamics reduced-order models},
   volume = {0},
   year = {2009},
}

@article{bruls2007,
   author = {Olivier Br\"{u}ls and Pierre Duysinx and Jean Claude Golinval},
   doi = {10.1002/nme.1795},
   issn = {00295981},
   issue = {5},
   journal = {International Journal for Numerical Methods in Engineering},
   month = {1},
   pages = {948-977},
   title = {The global modal parameterization for non-linear model-order reduction in flexible multibody dynamics},
   volume = {69},
   year = {2007},
}

@article{kim2007,
   author = {Kim TS, Kim YY},
   journal = {Computers and Structures},
   title = {MAC-based mode-tracking in structural topology optimization},
   year = {2000},
}

@article{Stephen2009,
   author = {N. G. Stephen},
   doi = {10.1115/1.3147130},
   issn = {15288927},
   issue = {5},
   journal = {Journal of Vibration and Acoustics},
   pages = {0545011-0545015},
   publisher = {American Society of Mechanical Engineers(ASME)},
   title = {On veering of eigenvalue loci},
   volume = {131},
   year = {2009},
}

@article{Heirman2011,
   author = {Gert H.K. Heirman and Frank Naets and Wim Desmet},
   doi = {10.1002/nme.2971},
   issn = {00295981},
   issue = {3},
   journal = {International Journal for Numerical Methods in Engineering},
   month = {1},
   pages = {330-354},
   title = {A system-level model reduction technique for the efficient simulation of flexible multibody systems},
   volume = {85},
   year = {2011},
}

@inproceedings{bonisoli2013crossing,
  title={Crossing and veering phenomena in crank mechanism dynamics},
  author={Bonisoli, Elvio and Marcuccio, Gabriele and Rosso, Carlo},
  booktitle={Topics in Model Validation and Uncertainty Quantification, Volume 5: Proceedings of the 31st IMAC, A Conference on Structural Dynamics, 2013},
  pages={175--187},
  year={2013},
  organization={Springer}
}

@article{Giannini2016,
   author = {O. Giannini and A. Sestieri},
   doi = {10.1016/j.ymssp.2015.11.012},
   issn = {10961216},
   journal = {Mechanical Systems and Signal Processing},
   month = {5},
   pages = {846-864},
   publisher = {Academic Press},
   title = {Experimental characterization of veering crossing and lock-in in simple mechanical systems},
   volume = {72-73},
   year = {2016},
}

@article{allemang2003modal,
  title={The modal assurance criterion--twenty years of use and abuse},
  author={Allemang, Randall J},
  journal={Sound and vibration},
  volume={37},
  number={8},
  pages={14--23},
  year={2003}
}

@inproceedings{Amsallem2011,
   author = {David Amsallem and Charbel Farhat},
   doi = {10.1137/100813051},
   issn = {10648275},
   issue = {5},
   journal = {SIAM Journal on Scientific Computing},
   pages = {2169-2198},
   publisher = {Society for Industrial and Applied Mathematics Publications},
   title = {An online method for interpolating linear parametric reduced-order models},
   volume = {33},
   year = {2011},
}

@article{Saccani2024,
   author = {Alexander Saccani and Paolo Tiso},
   month = {11},
   title = {Accelerating Construction of Non-Intrusive Nonlinear Structural Dynamics Reduced Order Models through Hyperreduction},
   url = {http://arxiv.org/abs/2411.14262},
   year = {2024},
}

@article{Hollkamp2005,
   author = {Joseph J. Hollkamp and Robert W. Gordon and S. Michael Spottswood},
   doi = {10.1016/j.jsv.2004.08.036},
   issn = {0022460X},
   issue = {3-5},
   journal = {Journal of Sound and Vibration},
   month = {6},
   pages = {1145-1163},
   publisher = {Academic Press},
   title = {Nonlinear modal models for sonic fatigue response prediction: A comparison of methods},
   volume = {284},
   year = {2005},
}

@misc{Radu2003,
   author = {Adrian George Radu and Marc P Mignolet and Adrian G Radu and Xiaowei Gao},
   title = {VALIDATION OF REDUCED ORDER MODELING FOR THE PREDICTION OF THE RESPONSE AND FATIGUE LIFE OF PANELS SUBJECTED TO THERMO-ACOUSTIC EFFECTS},
   url = {https://www.researchgate.net/publication/262933555},
   year = {2003},
}

@article{mignolet2013review,
	title={A review of indirect/non-intrusive reduced order modeling of nonlinear geometric structures},
	author={Mignolet, Marc P and Przekop, Adam and Rizzi, Stephen A and Spottswood, S Michael},
	journal={Journal of Sound and Vibration},
	volume={332},
	number={10},
	pages={2437--2460},
	year={2013},
	publisher={Elsevier}
}

@BOOK{crisfield1991,
  TITLE = {Nonlinear Finite Element Analysis of Solids and Structures},
  AUTHOR = {M.A. Crisfield},
  YEAR = {1991},
  PUBLISHER = {Wiley},
}

@incollection{tiso20214,
  title={Modal methods for reduced order modeling},
  author={Tiso, Paolo and Karamooz Mahdiabadi, Morteza and Marconi, Jacopo},
  booktitle={Model Order Reduction: Volume 1: System-and Data-Driven Methods and Algorithms},
  pages={97--138},
  year={2021},
  publisher={De Gruyter}
}

@article{Kim2013,
   author = {Kwangkeun Kim and Adrian G. Radu and X. Q. Wang and Marc P. Mignolet},
   doi = {10.1016/j.ijnonlinmec.2012.07.008},
   issn = {00207462},
   journal = {International Journal of Non-Linear Mechanics},
   pages = {100-110},
   publisher = {Elsevier Ltd},
   title = {Nonlinear reduced order modeling of isotropic and functionally graded plates},
   volume = {49},
   year = {2013},
}

@article{Idelsohn1985,
   author = {Sergio R. Idelsohn and Alberto Cardona},
   doi = {10.1016/0045-7825(85)90125-2},
   issn = {00457825},
   issue = {3},
   journal = {Computer Methods in Applied Mechanics and Engineering},
   pages = {253-279},
   title = {A reduction method for nonlinear structural dynamic analysis},
   volume = {49},
   year = {1985},
}

@article{Marconi2020,
   author = {Jacopo Marconi and Paolo Tiso and Francesco Braghin},
   doi = {10.1016/j.cma.2019.112785},
   issn = {00457825},
   journal = {Computer Methods in Applied Mechanics and Engineering},
   month = {3},
   publisher = {Elsevier B.V.},
   title = {A nonlinear reduced order model with parametrized shape defects},
   volume = {360},
   year = {2020},
}

@article{Marconi2021,
   author = {J. Marconi and P. Tiso and D. E. Quadrelli and F. Braghin},
   doi = {10.1007/s11071-021-06496-y},
   issn = {1573269X},
   issue = {4},
   journal = {Nonlinear Dynamics},
   month = {6},
   pages = {3039-3063},
   publisher = {Springer Science and Business Media B.V.},
   title = {A higher-order parametric nonlinear reduced-order model for imperfect structures using Neumann expansion},
   volume = {104},
   year = {2021},
}

@article{Wu2016,
   author = {Long Wu and Paolo Tiso},
   doi = {10.1007/s11044-015-9476-5},
   issn = {1573272X},
   issue = {4},
   journal = {Multibody System Dynamics},
   month = {4},
   pages = {405-425},
   publisher = {Springer Netherlands},
   title = {Nonlinear model order reduction for flexible multibody dynamics: a modal derivatives approach},
   volume = {36},
   year = {2016},
}

@article{Sombroek2018,
   author = {C. S.M. Sombroek and P. Tiso and L. Renson and G. Kerschen},
   doi = {10.1016/j.compstruc.2017.08.016},
   issn = {00457949},
   journal = {Computers and Structures},
   month = {1},
   pages = {34-46},
   publisher = {Elsevier Ltd},
   title = {Numerical computation of nonlinear normal modes in a modal derivative subspace},
   volume = {195},
   year = {2018},
}

@article{Morteza,
   author = {Morteza Karamooz Mahdiabadi and Paolo Tiso and Antoine Brandt and Daniel Jean Rixen},
   doi = {10.1016/j.ymssp.2020.107126},
   issn = {10961216},
   journal = {Mechanical Systems and Signal Processing},
   month = {1},
   publisher = {Academic Press},
   title = {A non-intrusive model-order reduction of geometrically nonlinear structural dynamics using modal derivatives},
   volume = {147},
   year = {2021},
}

@article{idelsohn1985load,
  title={A load-dependent basis for reduced nonlinear structural dynamics},
  author={Idelsohn, Sergio R and Cardona, Alberto},
  journal={Computers \& Structures},
  volume={20},
  number={1-3},
  pages={203--210},
  year={1985},
  publisher={Elsevier}
}

@article{Andersson2023,
   author = {Linus Andersson and Peter Persson and Kent Persson},
   doi = {10.1016/j.ymssp.2023.110143},
   issn = {10961216},
   journal = {Mechanical Systems and Signal Processing},
   month = {5},
   publisher = {Academic Press},
   title = {Efficient nonlinear reduced order modeling for dynamic analysis of flat structures},
   volume = {191},
   year = {2023},
}

@article{Jain2017,
   author = {Shobhit Jain and Paolo Tiso and Johannes B. Rutzmoser and Daniel J. Rixen},
   doi = {10.1016/j.compstruc.2017.04.005},
   issn = {00457949},
   journal = {Computers and Structures},
   month = {8},
   pages = {80-94},
   publisher = {Elsevier Ltd},
   title = {A quadratic manifold for model order reduction of nonlinear structural dynamics},
   volume = {188},
   year = {2017},
}

@book{geradin2015mechanical,
  title={Mechanical vibrations: theory and application to structural dynamics},
  author={G{\'e}radin, Michel and Rixen, Daniel J},
  year={2015},
  publisher={John Wiley \& Sons}
}

@article{Farhat2015,
   author = {Charbel Farhat and Todd Chapman and Philip Avery},
   doi = {10.1002/nme.4820},
   issn = {10970207},
   issue = {5},
   journal = {International Journal for Numerical Methods in Engineering},
   month = {5},
   pages = {1077-1110},
   publisher = {John Wiley and Sons Ltd},
   title = {Structure-preserving, stability, and accuracy properties of the energy-conserving sampling and weighting method for the hyper reduction of nonlinear finite element dynamic models},
   volume = {102},
   year = {2015},
}

@article{Rutzmoser2017,
   author = {J. B. Rutzmoser and D. J. Rixen},
   doi = {10.1016/j.cma.2017.06.009},
   issn = {00457825},
   journal = {Computer Methods in Applied Mechanics and Engineering},
   month = {10},
   pages = {330-349},
   publisher = {Elsevier B.V.},
   title = {A lean and efficient snapshot generation technique for the Hyper-Reduction of nonlinear structural dynamics},
   volume = {325},
   year = {2017},
}

@article{Jain2018,
   author = {Shobhit Jain and Paolo Tiso},
   doi = {10.1115/1.4040021},
   issn = {15551423},
   issue = {7},
   journal = {Journal of Computational and Nonlinear Dynamics},
   month = {7},
   publisher = {American Society of Mechanical Engineers (ASME)},
   title = {Simulation-free hyper-reduction for geometrically nonlinear structural dynamics: A quadratic manifold lifting approach},
   volume = {13},
   year = {2018},
}

@phdthesis{rutzmoserThesis,
	author = {Rutzmoser, Johannes},
	title = {Model Order Reduction for Nonlinear Structural Dynamics},
	year = {2018},
	school = {Technische Universit\"at M\"unchen},
	pages = {257},
	language = {en}
}

@inproceedings{McEwan2001,
   author = {M. I. McEwan and J. R. Wright and J. E. Cooper and A. Y.T. Leung},
   doi = {10.2514/6.2001-1595},
   journal = {19th AIAA Applied Aerodynamics Conference},
   publisher = {American Institute of Aeronautics and Astronautics Inc.},
   title = {A finite element/modal technique for nonlinear plate and stiffened panel response prediction},
   year = {2001},
}

@article{Muravyov2003,
   author = {Alexander A. Muravyov and Stephen A. Rizzi},
   doi = {10.1016/S0045-7949(03)00145-7},
   issn = {00457949},
   issue = {15},
   journal = {Computers and Structures},
   month = {7},
   pages = {1513-1523},
   title = {Determination of nonlinear stiffness with application to random vibration of geometrically nonlinear structures},
   volume = {81},
   year = {2003},
}

@article{Perez2014,
   author = {Ricardo Perez and X. Q. Wang and Marc P. Mignolet},
   doi = {10.1115/1.4026155},
   issn = {15551415},
   issue = {3},
   journal = {Journal of Computational and Nonlinear Dynamics},
   month = {7},
   title = {Nonintrusive structural dynamic reduced order modeling for large deformations: Enhancements for complex structures},
   volume = {9},
   year = {2014},
}

@inproceedings{Gordon2011,
   author = {Robert W. Gordon and Joseph J. Hollkamp},
   doi = {10.2514/6.2011-2081},
   isbn = {9781600869518},
   issn = {02734508},
   journal = {Collection of Technical Papers - AIAA/ASME/ASCE/AHS/ASC Structures, Structural Dynamics and Materials Conference},
   title = {Reduced-order models for acoustic response prediction of a curved panel},
   year = {2011},
}

@inproceedings{Spottswood2008,
   author = {Stephen Spottswood and Joseph Hollkamp and Thomas Eason},
   doi = {10.2514/6.2008-2235},
   month = {4},
   publisher = {American Institute of Aeronautics and Astronautics (AIAA)},
   title = {On the Use of Reduced-Order Models for a Shallow Curved Beam Under Combined Loading},
   year = {2008},
}

@inproceedings{tiso2011optimal,
  title={Optimal second order reduction basis selection for nonlinear transient analysis},
  author={Tiso, Paolo},
  booktitle={Modal Analysis Topics, Volume 3: Proceedings of the 29th IMAC, A Conference on Structural Dynamics, 2011},
  pages={27--39},
  year={2011},
  organization={Springer}
}

@article{helton2003latin,
  title={Latin hypercube sampling and the propagation of uncertainty in analyses of complex systems},
  author={Helton, Jon C and Davis, Freddie Joe},
  journal={Reliability Engineering \& System Safety},
  volume={81},
  number={1},
  pages={23--69},
  year={2003},
  publisher={Elsevier}
}

@inproceedings{przekop2006nonlinear,
  title={Nonlinear acoustic response of an aircraft fuselage sidewall structure by a reduced-order analysis},
  author={Przekop, Adam and Rizzi, Stephen A and Groen, David S},
  booktitle={Ninth International Conference on Recent Advances in Structural Dynamics},
  number={Paper-135},
  year={2006}
}

@article{Przekop2007,
   author = {Adam Przekop and Stephen A. Rizzi},
   doi = {10.2514/1.26351},
   issn = {00011452},
   issue = {10},
   journal = {AIAA Journal},
   month = {10},
   pages = {2510-2519},
   title = {Dynamic snap-through of thin-walled structures by a reduced-order method},
   volume = {45},
   year = {2007},
}

@article{Park2023,
   author = {Kyusic Park and Matthew S. Allen},
   doi = {10.1016/j.ymssp.2022.109720},
   issn = {10961216},
   journal = {Mechanical Systems and Signal Processing},
   month = {2},
   publisher = {Academic Press},
   title = {A Gaussian process regression reduced order model for geometrically nonlinear structures},
   volume = {184},
   year = {2023},
}

@article{Hill2016,
   author = {T. L. Hill and P. L. Green and A. Cammarano and S. A. Neild},
   doi = {10.1016/j.jsv.2015.09.007},
   issn = {10958568},
   journal = {Journal of Sound and Vibration},
   month = {1},
   pages = {156-170},
   publisher = {Academic Press},
   title = {Fast Bayesian identification of a class of elastic weakly nonlinear systems using backbone curves},
   volume = {360},
   year = {2016},
}

@article{Song2018,
   author = {Mingming Song and Ludovic Renson and Jean Philippe No\"{e}l and Babak Moaveni and Gaetan Kerschen},
   doi = {10.1002/stc.2258},
   issn = {15452263},
   issue = {12},
   journal = {Structural Control and Health Monitoring},
   month = {12},
   publisher = {John Wiley and Sons Ltd},
   title = {Bayesian model updating of nonlinear systems using nonlinear normal modes},
   volume = {25},
   year = {2018},
}

@article{Dubois,
   author = {A. Paul-Dubois-Taine and D. Amsallem},
   doi = {10.1002/nme.4759},
   issn = {10970207},
   issue = {5},
   journal = {International Journal for Numerical Methods in Engineering},
   month = {5},
   pages = {1262-1292},
   publisher = {John Wiley and Sons Ltd},
   title = {An adaptive and efficient greedy procedure for the optimal training of parametric reduced-order models},
   volume = {102},
   year = {2015},
}

@book{wright2003radial,
  title={Radial basis function interpolation: numerical and analytical developments},
  author={Wright, Grady Barrett},
  year={2003},
  publisher={University of Colorado at Boulder}
}

@software{yafec,
  author       = {S. Jain and J. Marconi and P. Tiso},
  title        = {YetAnotherFeCode},
  year         = {2021},
  version      = {v 1.1.1}
}

\end{document}